\documentclass[12pt]{article}
\usepackage{amsmath}
\usepackage{amssymb}
\usepackage{amsthm}
\usepackage[totalwidth=17cm, totalheight=25cm]{geometry}
\usepackage{graphicx}
\usepackage{picins}
\usepackage{tabularx}
	\newcolumntype{C}[1]{>{\centering\arraybackslash}m{#1}} 
	\newcolumntype{R}[1]{>{\raggedleft\arraybackslash}m{#1}} 

\newtheoremstyle{boldplain}
{9pt}
{9pt}
{\itshape}
{}
{\bfseries}
{.}
{.5em}
{\thmname{#1}\thmnumber{ #2}\thmnote{ (#3)}}%

\newtheoremstyle{bolddefinition}
{9pt}
{9pt}
{}
{}
{\bfseries}
{.}
{.5em}
{\thmname{#1}\thmnumber{ #2}\thmnote{ (#3)}}%

\theoremstyle{boldplain}

\newtheorem{conj}[equation]{Conjecture}
\newtheorem{cor}[equation]{Corollary}
\newtheorem{lem}[equation]{Lemma}
\newtheorem{prop}[equation]{Proposition}
\newtheorem{ques}[equation]{Question}
\newtheorem{sublem}[equation]{Sublemma}
\newtheorem{thm}[equation]{Theorem}

\theoremstyle{bolddefinition}

\newtheorem{rem}[equation]{Remark}

\newfont{\bigbf}{cmbx10 scaled\magstep1}
\normalbaselines
\numberwithin{equation}{section}

\def\no{\noindent}

\def\R{{\mathbb R}}

\def\Z{{\mathbb Z}}

\def\al{\alpha}
\def\ga{\gamma}
\def\Ga{\Gamma}
\def\de{\delta}
\def\De{\Delta}
\def\eps{\epsilon}

\def\si{\sigma}
\def\Si{\Sigma}

\def\3{\ss}

\def\acts{\curvearrowright}
\def\embed{\hookrightarrow}
\def\half{\frac{1}{2}}
\def\third{\frac{1}{3}}

\def\ol{\overline}
\def\ora{\overrightarrow}
\def\pihalf{\frac{\pi}{2}}
\def\pithird{\frac{\pi}{3}}
\def\2pithird{\frac{2\pi}{3}}
\def\piforth{\frac{\pi}{4}}
\def\pisixth{\frac{\pi}{6}}
\def\quart{\frac{1}{4}}

\def\phm{\phantom{-}}

\def\diam{\mathop{\hbox{diam}}}
\def\rad{\mathop{\hbox{rad}}}

\title{The Center Conjecture for spherical buildings \\
of types $F_4$ and $E_6$} 
\author{Bernhard Leeb, Carlos Ramos-Cuevas\footnote{
b.l@lmu.de, cramos@mathematik.uni-muenchen.de}}
\date{February 4, 2011}

\begin{document}

\maketitle
\centerline{\it Dedicated to Jos\'e Mar\'ia
 Montesinos on the occasion of his 65th birthday}

\begin{abstract}
\noindent
We prove that a convex subcomplex 
of a spherical building of type $F_4$ or $E_6$
is a subbuilding 
or the automorphisms of the subcomplex fix a point on it. 
Our approach is differential-geometric 
and based on the theory of metric spaces with curvature bounded above.
We use these techniques also 
to give another proof of the same result 
for the spherical buildings of classical type. 
\end{abstract}

\section{Introduction}

A subset in a CAT(1) space,
i.e.\ in a space with curvature $\leq1$ in the comparison sense, 
is called convex 
if it contains with any two points of distance $<\pi$ 
the unique minimizing geodesic segment connecting them. 

In a Euclidean unit sphere 
there are no proper convex subsets beyond a certain threshold:
A convex subset is either contained in a convex metric ball,
that is, a ball of radius $\leq\frac{\pi}{2}$
or it fills out the entire sphere. 
To put it more intrinsically, 
the convex subset is either contained in a convex metric ball
centered around one of its points 
or it is a geodesic subsphere.
Thus its intrinsic circumradius is $\leq\pihalf$ or $=\pi$.

Spherical buildings are a very special kind of CAT(1) spaces. 
Their geometry is rigidified by the property 
that they contain ``plenty of apartments'',
i.e.\ top-dimensional convex subsets isometric to a unit sphere.
A metric ball of radius $<\pi$ 
in a spherical building is convex if and only if 
it has radius $\leq\frac{\pi}{2}$. 
It is natural to ask whether the ``circumradius gap phenomenon'' 
for convex subsets of spheres 
holds more generally in spherical buildings, 
compare \cite[Question 1.5]{invconv}:

\begin{ques}
\label{gapconj}
Suppose that $B$ is a spherical building 
and that $C\subseteq B$ is a convex subset. 
Is it true that $C$ is either a subbuilding 
or it is contained in a convex metric ball
centered in $C$? 
\end{ques}

It is easy to see that the answer is yes, if $dim(C)\leq1$.
A one-dimensional convex subset is either a subbuilding or a tree.
In the latter case, 
it contains a unique circumcenter. 

Regarding isometric group actions on spherical buildings, 
one can ask the following weaker version of Question~\ref{gapconj}, 
see \cite{BalserLytchak_centers}: 

\begin{ques}
\label{fixedpt}
Suppose that $B$ is a spherical building 
and that $C\subseteq B$ is a convex subset. 
Is it true that $C$ is a subbuilding 
or the action $Isom(C)\acts C$ has a fixed point? 
\end{ques}

A positive answer to Question~\ref{gapconj}
implies a positive answer to Question~\ref{fixedpt}. 
 
Question~\ref{fixedpt} has been answered positively in 
\cite[Thm.\ 1.1] {BalserLytchak_centers}
when $dim(C)\leq2$. 

In higher dimensions, 
both questions seem to become considerably more approachable
when one restricts to convex subsets 
which are {\em subcomplexes} 
with respect to the natural polyhedral structure 
of the spherical building. 
Question~\ref{fixedpt} then becomes 
a geometric version of Tits' Center Conjecture, 
compare 
\cite{MuehlherrTits} and \cite[Conjecture 2.8]{Serre}: 

\begin{conj}[Center Conjecture]
\label{conj:tcc}
Suppose that $B$ is a spherical building 
and that $K\subseteq B$ is a convex subcomplex. 
Then $K$ is a subbuilding 
or the action $Stab_{Aut(B)}(K)\acts K$ 
of the automorphisms of $B$ preserving $K$ has a fixed point. 
\end{conj}
The automorphisms of a spherical building are the isometries 
which preserve its combinatorial (i.e.\ polyhedral) structure. 

A positive answer to Question~\ref{fixedpt} 
implies a positive answer to Conjecture~\ref{conj:tcc}.

The Center Conjecture was first proposed by Tits in the 1950s. 
Part of his motivation came from algebraic group theory, 
cf.\ \cite[Lemma 1.2]{Tits_cc}. 
A special case of the conjecture is also considered in 
Geometric Invariant Theory
related to the discussion of instability, see
\cite{Mumford},
and this special case had been proven by 
Rousseau \cite{Rousseau} and Kempf \cite{Kempf}. 

Conjecture~\ref{conj:tcc} easily reduces to the irreducible case. 
It has been proven for irreducible buildings of types $A_n$, $B_n$ and $D_n$ 
in \cite{MuehlherrTits}. 
The $F_4$-case has been announced by Parker and Tent
in a talk at Oberwolfach 
 in January 2008 \cite{ParkerTent_conv}.
To our knowledge no written account of the argument is available.
These approaches are incidence-geometric. 
 
The main result of this paper 
is the proof of the Center Conjecture~\ref{conj:tcc} 
for spherical buildings of types $F_4$ and $E_6$,
see Theorems~\ref{thm:tccf4} and \ref{thm:tcce6}. 
Our methods are differential-geometric 
and based on the theory of metric spaces with curvature bounded above 
and in particular on the geometric approach to buildings 
from the perspective of comparison geometry 
as in \cite{qirigid}. 
The arguments rely on the specific features 
of $F_4$- and $E_6$-geometry. 
We use these techniques also 
to give another proof of the Center Conjecture 
for the spherical buildings of classical type. 
Our arguments actually yield a slightly stronger intrinsic
version of the Center Conjecture where we admit all
intrinsic automorphisms of the subcomplex $K$ instead of only those
which extend to automorphisms of the ambient building $B$,
see Theorem~\ref{thm:generalversion}.

The approach in this paper has been carried further by the second author 
in \cite{E8,diss} 
where he proves Conjecture~\ref{conj:tcc} 
for spherical buildings of types $E_7$ and $E_8$. 
It follows that the Center Conjecture holds 
for all spherical buildings without factors of type $H_4$,
and in particular for all thick spherical buildings. 

\medskip
{\em Acknowledgements.}
The first author thanks C.\ Parker and K.\ Tent 
for reviving his interest in the Center Conjecture. 
We are grateful to the CONACYT (Mexico), 
the Deutscher Akademischer Austauschdienst, 
the Deutsche Forschungsgemeinschaft
and the Hausdorff Research Institute for Mathematics 
for financial support,
and to the Hausdorff Institute 
also for its great hospitality during our stay in fall 2008.

\tableofcontents

\section{Some geometric properties of spherical buildings}
\label{sec:prelim}

\subsection{CAT(1) spaces}
\label{sec:cat1} 

A complete metric space $X$ is called a {\em CAT(1) space}, 
if any two points with distance $<\pi$ 
are connected by a minimizing geodesic segment 
and if geodesic triangles with perimeter $<2\pi$ 
are not thicker than the corresponding comparison triangles 
in the unit sphere $S^2(1)$ with Gau\3 curvature $\equiv1$. 
We refer to 
\cite[ch.\ 2.1-2]{qirigid} and \cite[ch.\ 2.1-3]{BridsonHaefliger} 
for basic information on CAT(1) spaces. 

We denote by $B_r(x)$ the open metric ball 
of radius $r$ centered at $x$, 
and by $xy$ a minimizing geodesic segment 
with endpoints $x$ and $y$. 

The {\em link} or {\em space of directions} $\Si_xX$
of $X$ at a point $x$ 
equipped with the angle metric 
is again a CAT(1) space. 
It can be thought of as an analogue 
of the unit tangent sphere of a Riemannian manifold. 
If $x\neq y$, 
we denote by $\ora{xy}\in\Si_xX$ 
the direction of the segment $xy$ at $x$.

\subsubsection{Convexity}

One calls a subset $C\subseteq X$ 
{\em $\pi$-convex} or simply {\em convex}, 
if with any two points $x,x'\in C$ of distance $<\pi$ 
the unique minimizing geodesic segment $xx'$ 
is contained in $C$. 
Closed convex subsets of CAT(1) spaces 
are CAT(1) spaces themselves. 
Metric balls with radius $\leq\pihalf$ in CAT(1) spaces are convex. 
The {\em closed convex hull} $CH(A)$ of a subset $A\subseteq X$ 
is the smallest closed convex subset of $X$ containing $A$. 
We will denote the closed convex hull of finitely many points 
$a_1,\dots,a_m$ by $CH(a_1,\dots,a_m)$. 

\subsubsection{Circumradius and circumcenters}
\label{sec:rad}

For a subset $A\subseteq X$ and a point $x\in X$ 
we denote by $\rad(A,x)$ 
the radius of the smallest closed metric ball around $x$ 
which contains $A$. 
We define the {\em circumradius} $\rad(A)=rad_X(A)$ of $A$ 
as the infimum of the function $\rad(A,\cdot)$ on $X$. 
A point where the infimum is attained 
is called a {\em circumcenter} of $A$.
If $A$ is convex, 
we call the infimum of $\rad(A,\cdot)$ {\em on $A$}
the {\em intrinsic} circumradius of $A$. 

If $\rad(A)<\frac{\pi}{2}$, 
then by standard comparison arguments 
$A$ has a unique circumcenter 
which must be contained in the closed convex hull $CH(A)$ of $A$. 
(Indeed, 
suppose that $(x_n)$ is a sequence of points in $X$ with 
$\rad(A,x_n)\searrow \rad(A)$,
and let $m_{ij}$ be the midpoint of $x_ix_j$. 
Then $\rad(A,m_{ij})\geq \rad(A)$ and the CAT(1) inequality 
imply that $(x_n)$ is a Cauchy sequence.
Its limit is a circumcenter of $A$, 
and it must be unique. 
The circumcenter must belong to $CH(A)$ 
because due to the CAT(1) inequality 
its nearest point projection to $CH(A)$ is also a circumcenter.)

If $\rad(A)=\pihalf$, 
then the set $Cent(A)$ of circumcenters of $A$ is closed and convex. 
(Its convexity follows from the CAT(1) inequality.)
Since $\rad(CH(A))=\pihalf$ and $Cent(A)=Cent(CH(A))$, 
the closed convex set $CH(A)\cap Cent(A)$ has diameter $\leq\pihalf$.

\subsection{Spherical Coxeter complexes}

We refer to 
\cite[ch.\ 4-5]{GroveBenson}, 
\cite[ch.\ V, VI.4]{Bourbaki}
and \cite[ch.\ 3.1, 3.3]{qirigid}
for more information. 

\subsubsection{General definitions and facts}
\label{sec:cox}

Let $S$ be the unit sphere 
in a finite dimensional Euclidean vector space $V$. 
The reflection at a hyperplane in $V$ through the origin 
induces an involutive isometry of $S$. 
One refers to such isometries briefly as {\em reflections}. 
If $W\subset Isom(S)$ is a finite subgroup generated by reflections, 
one calls the pair $(S,W)$ a {\em (spherical) Coxeter complex} 
and $W$ its {\em Weyl group}. 
(Note that we allow $W$ to have fixed points.)

The Weyl group $W$ induces a {\em polyhedral structure} on $S$. 
The fixed point sets of the reflections in $W$ 
are great spheres of codimension one, 
the {\em walls}.
There are finitely many walls 
and they divide $S$ into open convex subsets
whose closures are called {\em chambers}. 
If $W$ is nontrivial,
then the chambers are convex spherical polyhedra
because they are finite intersections of closed hemispheres. 
A {\em half-apartment} or {\em root} is a hemisphere bounded by a wall, 
a {\em singular sphere} is an intersection of walls, 
a {\em face} of $S$ 
is the intersection of a chamber with a singular sphere, 
a {\em panel} is a codimension one face, 
a {\em vertex} is a zero-dimensional face. 
Two faces are called {\em opposite} or {\em antipodal} 
if they are exchanged by the antipodal involution of $S$. 
The face {\em spanned} by a point 
is the face containing it as an interior point, 
equivalently,
the smallest face containing it. 
A point is called {\em regular} 
if it spans a chamber, 
and {\em singular} otherwise. 
A minimizing geodesic segment connecting two vertices 
is called {\em singular}
if it is contained in a singular 1-sphere.
A vertex is called of {\em root type} 
if the hemisphere centered at it is a root. 

Each chamber $\De$ is a fundamental domain 
for the action $W\acts S$ 
and $W$ is generated by the reflections at the codimension one faces of $\De$, 
that is, at the walls containing them. 
We call $\De_{mod}=\De_{mod}^{(S,W)}:=S/W$ 
the {\em model Weyl chamber}. 
Its isometry type determines $W$ up to conjugacy. 
The quotient map $\theta_S:S\to\De_{mod}$ is 1-Lipschitz 
and restricts to isometries on chambers. 
We call the image $\theta_S(x)$ of a point $x\in S$ its 
{\em $\theta_S$-type} or just its {\em type}. 

The {\em link} $\Si_xS$ of a point $x\in S$ 
is the unit tangent sphere of $S$ at $x$ 
in the sense of Riemannian geometry. 
It inherits from $S$ 
a natural structure as the spherical Coxeter complex 
$(\Si_xS,Stab_W(x))$ 
with Weyl group $Stab_W(x)$
and with model Weyl chamber 
$\De_{mod}^{(\Si_xS,Stab_W(x))}\cong\Si_x\De_{mod}^{(S,W)}$. 

More generally, 
let $\si\subset S$ be a face of codimension $\geq1$. 
Then for an interior point $x\in\si$, 
the link $\Si_xS$ splits as the spherical join 
$\Si_xS\cong\Si_x\si\circ\nu_x\si$
of the 
unit tangent sphere $\Si_x\si$ of $\si$ 
and the unit normal sphere $\nu_x\si$ of $\si$ in $S$. 
The unit normal sphere has dimension 
$\dim(S)-\dim(\si)-1$, 
and there is a natural isometric identification $\nu_x\si\cong Poles(\si)$
of $\nu_x\si$ with the sphere 
$Poles(\si):=\{p\in S:d(p,\cdot)|_{\si}\equiv\pihalf\}$ 
of {\em poles} of $\si$ in $S$. 
This provides one way to see, 
that one can consistently identify with each other 
the normal spheres $\nu_x\si$ for all interior points $x\in\si$ 
to obtain the {\em link $\Si_{\si}S$ of the face $\si$}. 
It inherits a natural structure as the spherical Coxeter complex 
$(\Si_{\si}S,Stab_W(\si))$ 
with Weyl group the stabilizer (fixator) of $\si$ in $W$ 
and with model Weyl chamber 
$\De_{mod}^{(\Si_{\si}S,Stab_W(\si))}\cong\Si_{\si}\De_{mod}^{(S,W)}$.

Let $s\subset S$ be a singular sphere. 
Then $s$ inherits a natural structure as a Coxeter complex as follows.
By a {\em reflection on $s$} we mean an involutive isometry of $s$
whose fixed point set is a codimension one subsphere. 
We define the {\em induced Weyl group} $W_s\subset Isom(s)$ on $s$ 
as the subgroup generated 
by those reflections on $s$ 
which are induced by isometries in $W$. 
The pair $(s,W_s)$ is a Coxeter complex 
and we refer to it as a {\em Coxeter subcomplex} of $(S,W)$. 
The Coxeter tesselation of $s$ is in general coarser
than its polyhedral structure inherited from $S$. 
Let us call the fixed point set of a reflection on $s$ and in $W_s$ 
an {\em $s$-wall}. 
Every face of codimension $\geq1$ in $s$ 
with respect to the (coarser) intrinsic polyhedral structure 
is contained in an $s$-wall. 
A codimension one face in $s$ 
with respect to the (finer) polyhedral structure induced from $S$ 
is contained in an $s$-wall 
if and only if both top-dimensional faces in $s$ adjacent to it 
(again with respect to the finer polyhedral structure) 
have the same type, i.e.\ the same $\theta_S$-image. 
\begin{rem}
Note that $W_s$ can be strictly smaller than 
the image of the natural homomorphism $Stab_W(s)\to Isom(s)$. 
An example for this phenomenon can be seen in 
the $E_7$-Coxeter complex: 
It contains a singular 1-sphere $s$ of type $13756137561$ 
(the first and last 1 to be identified).
The induced Weyl group $W_s$ is trivial, 
but the antipodal involution on $s$ 
is induced by isometries in $Stab_W(s)$. 
Here we use the labelling 
\hpic{\includegraphics[scale=0.4]{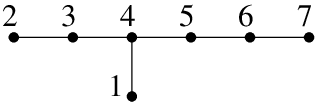}}
for the Dynkin diagram.
\end{rem}

We call a Coxeter complex {\em trivial} or a {\em sphere}
if its Weyl group is trivial. 
The Coxeter complex $(S,W)$ splits off a sphere factor 
if and only if $W$ has fixed points, 
equivalently, 
if and only if $\De_{mod}$ has diameter $\pi$. 
In this case the sphere $Fix(W)\subseteq S$ 
is canonically identified 
with the unique maximal sphere factor of $(S,W)$, 
its {\em spherical de Rham factor}.

We call a Coxeter complex {\em reducible} 
if it decomposes as the spherical join of Coxeter complexes. 
Join decompositions of a Coxeter complex 
correspond to join decompositions of its model Weyl chamber. 
For a Coxeter complex without spherical factor holds 
$\diam(\De_{mod})\leq\pihalf$. 
(If two $W$-orbits in $S$ have distance $>\pihalf$ 
then each of them is contained in an open hemisphere 
and has a center fixed by $W$.) 
It is irreducible if and only if 
$\diam(\De_{mod})<\pihalf$. 
The {\em de Rham decomposition} of $(S,W)$ 
is the unique maximal decomposition 
as the join of the spherical de Rham factor 
and some irreducible nontrivial Coxeter complexes. 

Suppose now that $W$ has no fixed points on $S$. 
Then $(S,W)$ has no spherical factor 
and $\De_{mod}$ is a spherical simplex. 
As remarked above, 
$\De_{mod}$ has diameter $\leq\pihalf$ 
with equality if and only if $(S,W)$ is reducible. 
The dihedral angle between any two panels of $\De_{mod}$ 
equals $\frac{\pi}{p}$ for some integer $p\geq2$. 
If $\De_{mod}$ has no one-dimensional join factor, 
then the only possible values for $p$ are 2,3,4 and 5. 
The geometry of $\De_{mod}$ can be encoded 
in a marked (or weighted) graph $\Ga$, 
the {\em Coxeter graph}, as follows. 
The vertices of $\Ga$ correspond to the panels of $\De_{mod}$. 
Two vertices are not connected if the corresponding dihedral angle 
equals $\pihalf$;
they are connected by an edge
if the angle equals $\pithird$,
and by an edge with label $p$ 
if the angle equals $\frac{\pi}{p}$ with $p\geq4$. 
If no edge labels $\neq4,6$ occur,
as it is the case for the Coxeter complexes coming from a
root system, 
then one often replaces the edges with label 4 by double edges
and the edges with label 6 by triple edges. 
The resulting graph with multiple edges 
is called the {\em Dynkin diagram}. 
The Coxeter graph determines $\De_{mod}$ up to isometry. 
Note that $\Ga$ is disconnected if and only if $(S,W)$ is reducible, 
equivalently, 
if $\De_{mod}$ decomposes as a spherical join. 
The classification of irreducible spherical Coxeter complexes
can be found in \cite[Thm.\ 5.3.1]{GroveBenson}. 
The irreducible Coxeter complexes of dimension $\geq2$ 
which occur for thick spherical buildings 
have the Dynkin diagrams 
$A_{n\geq3}$, $B_{n\geq3}$, $D_{n\geq4}$,  
$F_4$, $E_6$, $E_7$ and $E_8$, 
cf.\ \cite{Tits_w} and \cite[pp.\ 274]{Tits_bn}. 
For a face $\si\subset S$ of codimension $\geq1$, 
the Coxeter graph of its link $(\Si_{\si}S,Stab_W(\si))$ 
is obtained from the Coxeter graph of $(S,W)$ 
by deleting those vertices which correspond to the vertices of $\si$. 

When discussing a concrete Coxeter complex 
we will label the vertices of its Coxeter graph by some index set $I$. 
This induces also 
a labelling of the vertices of the Weyl chamber $\De_{mod}$ 
by assigning to a vertex $v$ of $\De_{mod}$ 
the label of the vertex of $\Ga$ corresponding 
to the panel opposite to $v$. 

An {\em automorphism} of the Coxeter complex $(S,W)$ 
is an isometry $\al$ of $S$ 
which preserves the tesselation into chambers, 
equivalently, which normalizes $W$,
i.e.\ $Aut((S,W))=N_{Isom(S)}(W)$. 
The isometries in $W$ are the {\em inner} automorphisms
and the automorphisms outside $W$ 
are the {\em outer} automorphisms.  
The {\em outer automorphism group}
$Out((S,W)):=Aut((S,W))/Inn((S,W))$ $=N_{Isom(S)}(W)/W$
is canonically identified with 
$Isom(\De_{mod})$
and with the automorphism group $Aut(\Ga)$ of the Coxeter graph. 

The antipodal involution of $S$ is always an automorphism of $(S,W)$. 
It induces the {\em canonical involution} $\iota$ of $\De_{mod}$.
For a chamber $\De\subset S$ 
there is a unique Weyl isometry $w\in W$ with $w\De=-\De$.
The composition $-w$ of $w$ with the antipodal involution of $S$ 
is an isometric involution of $\De$. 
It coincides with $\iota$ modulo the natural identification 
$\De\buildrel\theta_S\over\to\De_{mod}$. 
Note that $\iota$ is trivial if and only if 
any two opposite vertices in the Coxeter complex 
have equal type. 
Regarding the irreducible Coxeter complexes 
one has that 
$\iota=id_{\De_{mod}}$ for the Coxeter complexes of types 
$A_1$, $B_n$, $D_{2n}$, $F_4$, $E_7$ and $E_8$, 
and $\iota\neq id_{\De_{mod}}$ 
for the Coxeter complexes of types 
$A_{n\geq2}$, $D_{2n+1}$ and $E_6$.

\subsubsection{The Coxeter complex of type $F_4$}
\label{sec:f4geom}

Let $(S^3,W_{F_4})$ be the Coxeter complex of type $F_4$. 
We use the labelling 
\includegraphics[scale=0.6]{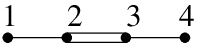}
for its Dynkin diagram $\Ga_{F_4}$. 
We collect here some geometric properties of $(S^3,W_{F_4})$ 
which will be needed in the proof of 
Theorem~\ref{thm:tccf4}
and which can be deduced from the information in 
\cite[ch.\ 5.3]{GroveBenson}.
To be more precise, the coordinate description of the
Coxeter complex will not be needed later but we give it here
in order to justify the other geometric properties.

We have 
$Out((S^3,W_{F_4}))\cong Isom(\De_{mod}^{F_4})
\cong Aut(\Ga_{F_4})\cong\Z_2$. 
The nontrivial involutive isometry of $\De_{mod}^{F_4}$
exchanges the vertices 
$1\leftrightarrow4$ and $2\leftrightarrow3$.
Hence the properties of $i$- and $(5-i)$-vertices 
are dual to each other. 

If we consider $(S^3,W_{F_4})$ embedded in 
$\R^4$ as the unit sphere,
we can describe the Weyl group as a group of 
isometries of $\R^4$ as follows.
The Weyl group $W_{F_4}$ is the finite group 
generated by the reflections at the hyperplanes orthogonal 
to the {\em fundamental root vectors}:
$$
r_1=-\half(1,1,1,1),\;\;\;\; 
r_2=e_1,\;\;\;\;
r_3=e_2-e_1,\;\;\;\;
r_2=e_3-e_2.
$$
The {\em fundamental Weyl chamber} $\De$ 
is given by the inequalities:
$$
x_1+\dots+x_4 \overset{\text{(1)}}{\leq} 0\;;\;\;\; \;\;
0\overset{\text{(2)}}{\leq} x_1
\overset{\text{(3)}}{\leq}x_2 
\overset{\text{(4)}}{\leq}x_3.
$$
We list vectors representing the vertices of $\De$:
\begin{center}
\begin{tabular}{p{3cm}p{6cm}}
1-vertex: $v_1$	&	$(\phm0,\phm0,\phm0,-1)$\\
2-vertex: $v_2$	&	$(\phm1,\phm1,\phm1,-3)$\\
3-vertex: $v_3$	&	$(\phm0,\phm1,\phm1,-2)$\\
4-vertex: $v_4$	&	$(\phm0,\phm0,\phm1,-1)$\\
\end{tabular}
\end{center}
All half-apartments of $(S^3,W_{F_4})$ are centered at a vertex. 
The vertices of types 1 and 4 are 
the vertices of {\em root type}. 
We list vectors representing these vertices:
\begin{center}
\begin{tabular}{C{2cm} C{8cm} }
1-vertices:	&	$\pm e_i$ for $i=1,\dots,4$;\phantom{..}	
				$\half(\pm e_1 \pm e_2 \pm e_3 \pm e_4)$\\ \\
4-vertices:	&	$\pm e_i \pm e_j$ for $1\leq i < j\leq 4$	
\end{tabular}
\end{center}
The vertices of root type are better separated from each other 
than the other types of vertices.
The possible mutual distances between 1-vertices
(4-vertices) 
are $0,\pithird,\pihalf,\2pithird$ and $\pi$. 
Any two pairs of 1-vertices with the same distance 
are equivalent modulo the action of the Weyl group.
Two 1-vertices with distance $\pithird$
are connected by a singular segment of type $121$, 
two 1-vertices with distance $\pihalf$
by a singular segment of type $141$, 
and two 1-vertices with distance $\2pithird$
by a singular segment of type $12121$.

\parpic(7cm,0pt)[l]{\includegraphics[scale=0.25]{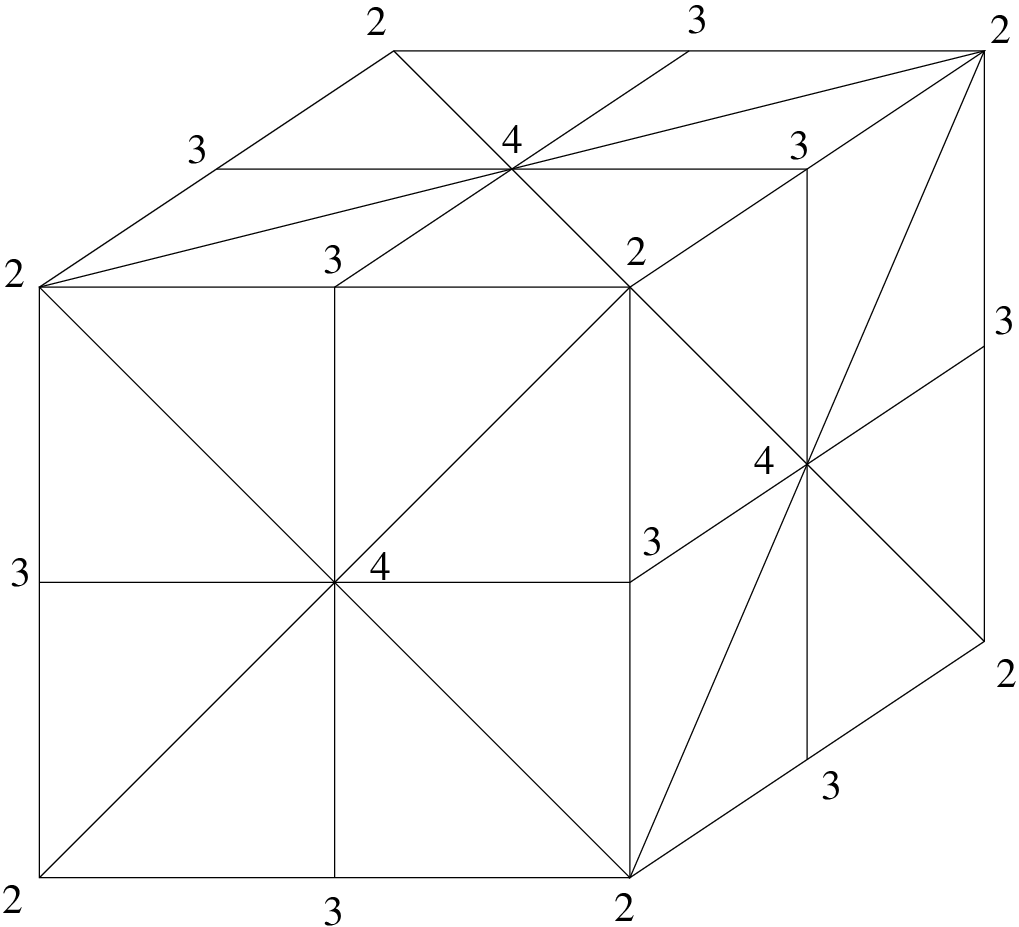}}
The link of a 1-vertex is a Coxeter complex 
$\Si_1$ of type $B_3$ 
with induced labelling 
\includegraphics[scale=0.6]{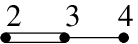}
for its Dynkin diagram.
Any two different 2-vertices in $\Si_1$ with distance $<\pihalf$
are connected by a singular segment of type $232$, 
and any two non-antipodal 2-vertices with distance $>\pihalf$
are connected by a singular segment of type $242$.

\parpic(7cm,3.7cm)[r]{\includegraphics[scale=0.2]{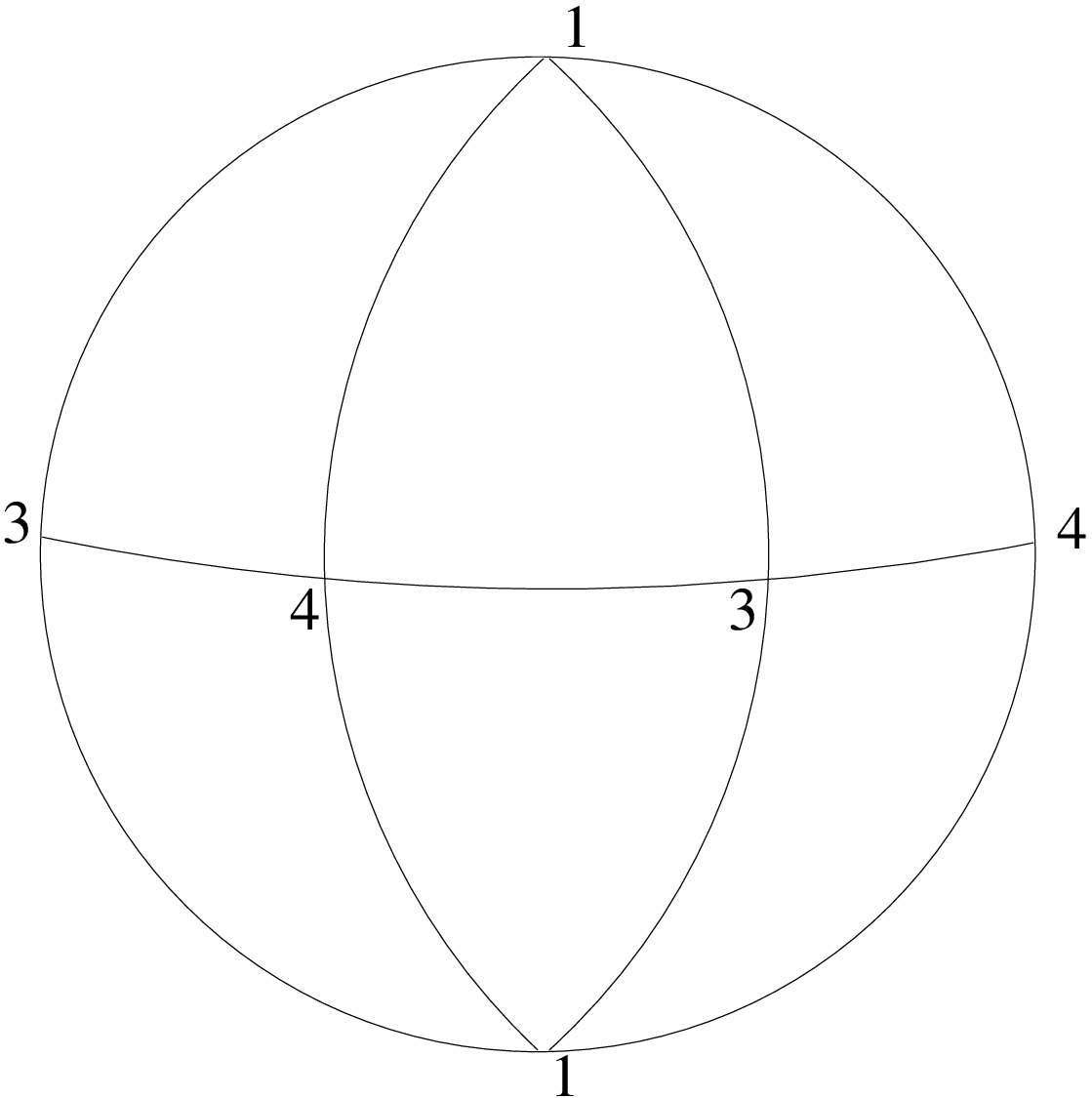}}
\vspace{1.5cm}
\picskip{3}
The link of a 2-vertex is a Coxeter complex 
$\Si_2$ of type $A_1\circ A_2$ 
with induced labelling 
\includegraphics[scale=0.6]{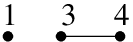}
for its Dynkin diagram. 
\vspace{0.8cm}

The following types of singular 1-spheres occur in $(S^3,W_{F_4})$:
\begin{center}
\includegraphics[scale=0.4]{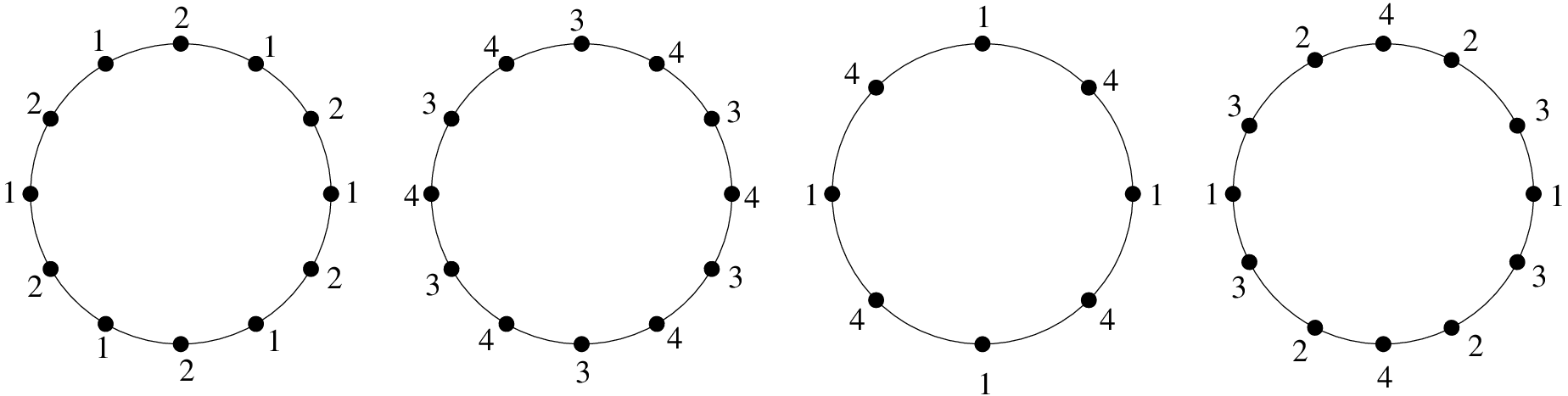}
\end{center}
The last type can be verified 
e.g.\ using the Dynkin diagrams of the links $\Si_i$,
compare section~\ref{sec:cox}, 
and the edge lengths of $\De$. 

The canonical involution 
$\iota:\De_{mod}^{F_4}\to\De_{mod}^{F_4}$ 
is trivial. 
Accordingly, the antipodes of $i$-vertices in the Coxeter complex 
are $i$-vertices. 
\begin{rem}
\label{rem:f4max}
Since 
$Out((S^3,W_{F_4}))\cong\Z_2$,
the normalizer $Aut((S^3,W_{F_4}))$ of $W_{F_4}$ in $Isom(S^3)$
is an index two extension of $W_{F_4}$. 
However,
it is not a reflection group,
because the nontrivial isometry of $\De$ fixes no vertex 
and therefore is not induced by a hyperplane reflection. 
\end{rem}

\subsubsection{The Coxeter complex of type $E_6$}
\label{sec:e6geom}

Let $(S^5,W_{E_6})$ be the Coxeter complex of type $E_6$. 
We use the labelling 
\hpic{\includegraphics[scale=0.4]{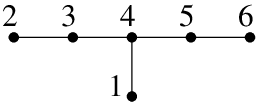}}
for its Dynkin diagram $\Ga_{E_6}$.
We collect here some geometric properties of $(S^5,W_{E_6})$ 
which will be needed in the proof of 
Theorem~\ref{thm:tcce6}
and which can be deduced from the information in 
\cite[ch.\ 5.3]{GroveBenson}. 
Again, the coordinate description of the
Coxeter complex will not be needed later but we give it here
in order to justify the other geometric properties.

We have 
$Out((S^5,W_{E_6}))\cong Isom(\De_{mod}^{E_6})
\cong Aut(\Ga_{E_6})\cong\Z_2$. 
The nontrivial involutive isometry of $\De_{mod}^{E_6}$
fixes the vertices 1 and 4 
and exchanges the vertices 
$2\leftrightarrow6$ and $3\leftrightarrow5$.

Our model for $(S^5,W_{E_6})$ is based on a model for the
Coxeter complex $(S^7,W_{E_8})$ of type $E_8$.
We consider $(S^7,W_{E_8})$ embedded in
$\R^8$ as the unit sphere
and will use the labelling
\hpic{\includegraphics[scale=0.4]{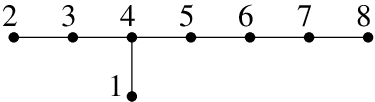}}
for its Dynkin diagram.
The root system of $E_8$ consists of the vectors
\begin{eqnarray*}
&\pm e_i\pm e_j
&\qquad\hbox{ for $1\leq i < j\leq 8$,}\\
&\half\sum\limits_{i=1}^8
\pm \eps_i e_i
&\qquad\hbox{ with $\eps_i=\pm1$ such that $\prod\limits_{i=1}^8\eps_i=-1$.}
\end{eqnarray*}
(The walls are the intersections of $S^7$ with the hyperplanes
perpendicular to a vector in the root system.)

The link $\Si_{\si}S^7$ of a type $78$ edge $\si$
is a Coxeter complex of type $E_6$.
The 8-vertices in $(S^7,W_{E_8})$ are the vertices of root type.
The $78$-edges have length $\pisixth$
and a pair of 8-vertices with distance $\pithird$ or $\2pithird$
lies on a singular circle of type $\dots8787\dots$.
Thus a model for $(S^5,W_{E_6})$ can be obtained
by choosing two $E_8$-root vectors $r,r'$ with angle $\2pithird$,
taking $S^5:=S^7\cap\langle r,r'\rangle^{\perp}$ 
and $W_{E_6}=Fix_{W_{E_8}}(\{r,r'\})$. 
The model in \cite{GroveBenson} uses 
$r=e_8-e_7$ and $r'=\half(1,1,1,1,1,1,1,-1)$.  
For us the choice 
$r=e_8-e_7$ and $r'=e_7-e_6$ 
is more convenient,
i.e.\ we realize $(S^5,W_{E_6})$ as the unit sphere in 
$\R^6\cong\{(x_1,\dots,x_8)\in\R^8\; |\; x_6=x_7=x_8\}$.

The $E_6$-root system then consists of the $E_8$-root vectors
perpendicular to $r$ and $r'$,
i.e.\ of 
\begin{equation}
\label{eq:e6roots}
\begin{array}{cccccc}
&\pm e_i\pm e_j
&\qquad\hbox{ for $1\leq i < j\leq 5$,}\hphantom{h that}&\\
&\half\sum\limits_{i=1}^8
\pm \eps_i e_i
&\qquad\hbox{ with $\eps_i=\pm1$ such that }
&\hbox{$\prod\limits_{i=1}^8\eps_i=-1$
 and $\eps_6=\eps_7=\eps_8$.}
\end{array}
\end{equation}
The Weyl group $W_{E_6}$
is the finite group of isometries
generated by the reflections at the hyperplanes orthogonal 
to the {\em fundamental root vectors}:
\begin{eqnarray*}
&r_1=\half(1,1,1,-1,-1,-1,-1,-1),\;\;\;\; r_i=e_i-e_{i-1} 
\text{ for } 2\leq i\leq 5; \\
&\text{ and }\;\;\; r_6=\half(1,1,1,1,-1,1,1,1).
\end{eqnarray*}
It contains as a proper subgroup the group $W'$ 
which permutes the first five coordinates 
and changes an even number of their signs. 

The {\em fundamental Weyl chamber} $\De$ 
is given by the inequalities:
\begin{equation*}
x_4+x_5+\dots+x_8 \overset{\text{(1)}}{\leq} x_1+x_2+x_3\;;\;\;\; \;\;
x_1\overset{\text{(2)}}{\leq} x_2
\overset{\text{(3)}}{\leq}\dots 
\overset{\text{(5)}}{\leq}x_5\;;\;\;\; \;\;
x_5\overset{\text{(6)}}{\leq}x_1+\dots+x_4+x_6+x_7+x_8.
\end{equation*}
We list vectors representing the vertices of $\De$:
\begin{center}
\begin{tabular}{p{3cm}p{6cm}}
1-vertex: $v_1$	&	$(\phm1,\phm1,\phm1,\phm1,\phm1,-1,-1,-1)$\\
2-vertex: $v_2$	&	$(-3,\phm3,\phm3,\phm3,\phm3,-1,-1,-1)$\\
3-vertex: $v_3$	&	$(\phm0,\phm0,\phm3,\phm3,\phm3,-1,-1,-1)$\\
4-vertex: $v_4$	&	$(\phm1,\phm1,\phm1,\phm3,\phm3,-1,-1,-1)$\\
5-vertex: $v_5$	&	$(\phm3,\phm3,\phm3,\phm3,\phm9,-1,-1,-1)$\\
6-vertex: $v_6$	&	$(\phm3,\phm3,\phm3,\phm3,\phm3,\phm1,\phm1,\phm1)$\\
\end{tabular}
\end{center}
All half-apartments in $(S^5,W_{E_6})$ are centered at a vertex. 
The {\em 1-vertices} are the vertices of {\em root type},
and they are represented by the root vectors (\ref{eq:e6roots}). 

The possible mutual distances between 1-vertices are 
$0$, $\pithird$, $\pihalf$, $\2pithird$ and $\pi$. 
The pairs of 1-vertices with the same distance 
are equivalent modulo the action of the Weyl group.
(This can be verified e.g.\ 
by considering the pairs containing $v_1$
which has a large stabilizer in $W'$.) 
A pair of 1-vertices with distance $\pithird$ ($\2pithird$) 
is connected by a type $141$ ($14141$) singular segment.
(Note that 
$v_4$ is the midpoint of $v_1$ and the 1-vertex
represented by $e_4+e_5$.) 
\parpic{\includegraphics[scale=0.4]{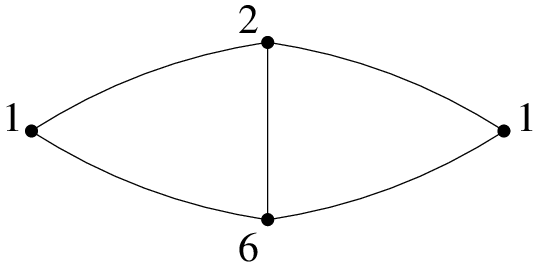}}
The segment connecting two 1-vertices with distance $\pihalf$
is not singular.
Its simplicial convex hull is a rhombus
whose other diagonal is a $26$-edge. 
(The midpoint of $v_2v_6$ equals the midpoint of $v_1$ 
and the 1-vertex represented by $\half(-1,1,\dots,1)$.)

An equilateral triangle with $141$-sides
is not a simplicial subcomplex. 
Its center is the midpoint of a $35$-edge 
perpendicular to the triangle. 
A square with $141$-sides
is not a simplicial subcomplex either.
Its center is the midpoint of a $26$-edge 
perpendicular to the square. 
\parpic{\includegraphics[scale=0.4]{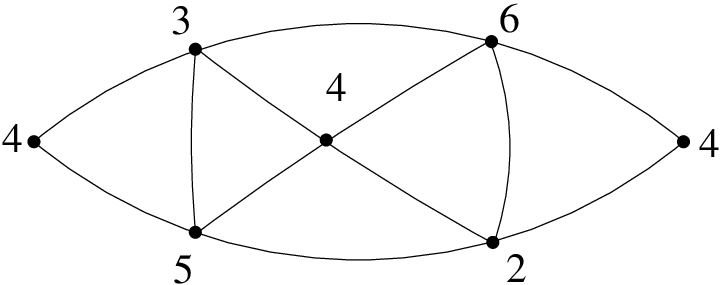}}
To verify these last facts, 
note that the link of a 1-vertex 
is a 4-dimensional Coxeter complex $\Si_1$ of type $A_5$ 
with induced labelling 
\includegraphics[scale=0.4]{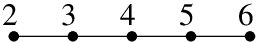}
for its Dynkin diagram. 
Any two distinct non-antipodal 4-vertices in $\Si_1$ 
have distance 
$\arccos(\pm\third)$, 
cf.\ section~\ref{app:coxeter}. 
The segment connecting them is not singular.
Its simplicial convex hull is a rhombus
whose other diagonal is a $35$-edge 
if the 4-vertices have distance $\arccos\third$,
and a $26$-edge if they have distance $\arccos(-\third)$.
Furthermore, 
the link of an edge of type $35$ ($26$)
in $(S^5,W_{E_6})$ 
is a Coxeter complex $\Si_{35}$ ($\Si_{26}$)
of type $A_2\circ A_1\circ A_1$ ($D_4$) 
with Dynkin diagram 
\includegraphics[scale=0.4]{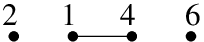}
(\hpic{\includegraphics[scale=0.3]{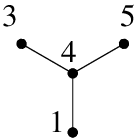}})
and contains a singular 1-sphere of type $\dots1414\dots$
with edge length $\pithird$ ($\piforth$). 

Next we list the {\em 2-vertices} 
modulo the action of $W'\subset W_{E_6}$,
more precisely, 
vectors representing them:
\begin{eqnarray*}
&(-3,\phm3,\phm3,\phm3,\phm3,-1,-1,-1)
,\;\;\;\; 
(\phm0,\phm0,\phm0,\phm0,\phm6,\phm2,\phm2,\phm2) \\
&(\phm0,\phm0,\phm0,\phm0,\phm0,-4,-4,-4)
\end{eqnarray*}
Since $-v_6$ is a 2-vertex,
the {\em 6-vertices} are just the antipodes of the 2-vertices.
We see that 
the canonical involution 
$\iota:\De_{mod}^{E_6}\to\De_{mod}^{E_6}$ 
is the nontrivial isometry. 
Accordingly, the antipodes of $i$-vertices in the Coxeter complex 
are $i$-vertices for $i=1,4$ 
and $(8-i)$-vertices for $i=2,3,5,6$. 

The possible mutual distances between 2-vertices
(6-vertices) 
are $0$, $\arccos\quart$ 
and $\2pithird$. 
Any two pairs of 2-vertices with the same distance 
are equivalent modulo the action of the Weyl group.
(This is obvious from considering the $W'$-orbits.)
Two 2-vertices (6-vertices) with distance $\arccos\quart$ 
are connected by a singular segment of type $232$ ($656$), 
and two 2-vertices (6-vertices) with distance $\2pithird$ 
by a singular segment of type $262$ ($626$). 
(One sees this by extending the edges $v_2v_3$ and $v_2v_6$ of $\De$.)

The possible distances between a 2-vertex and a 6-vertex 
are $\pithird$, $\arccos(-\quart)$ and $\pi$. 
If their distance is $\pithird$, 
they are connected by an edge of type $26$;
if their distance is $\arccos(-\quart)$, 
they are connected by a singular segment of type $216$. 
(This and the next list can be verified 
using the Dynkin diagrams of the links $\Si_i$ 
and the edge lengths of $\De$.)

The following types of singular 1-spheres occur in
$(S^5,W_{E_6})$: 
\begin{center}
\includegraphics[scale=0.4]{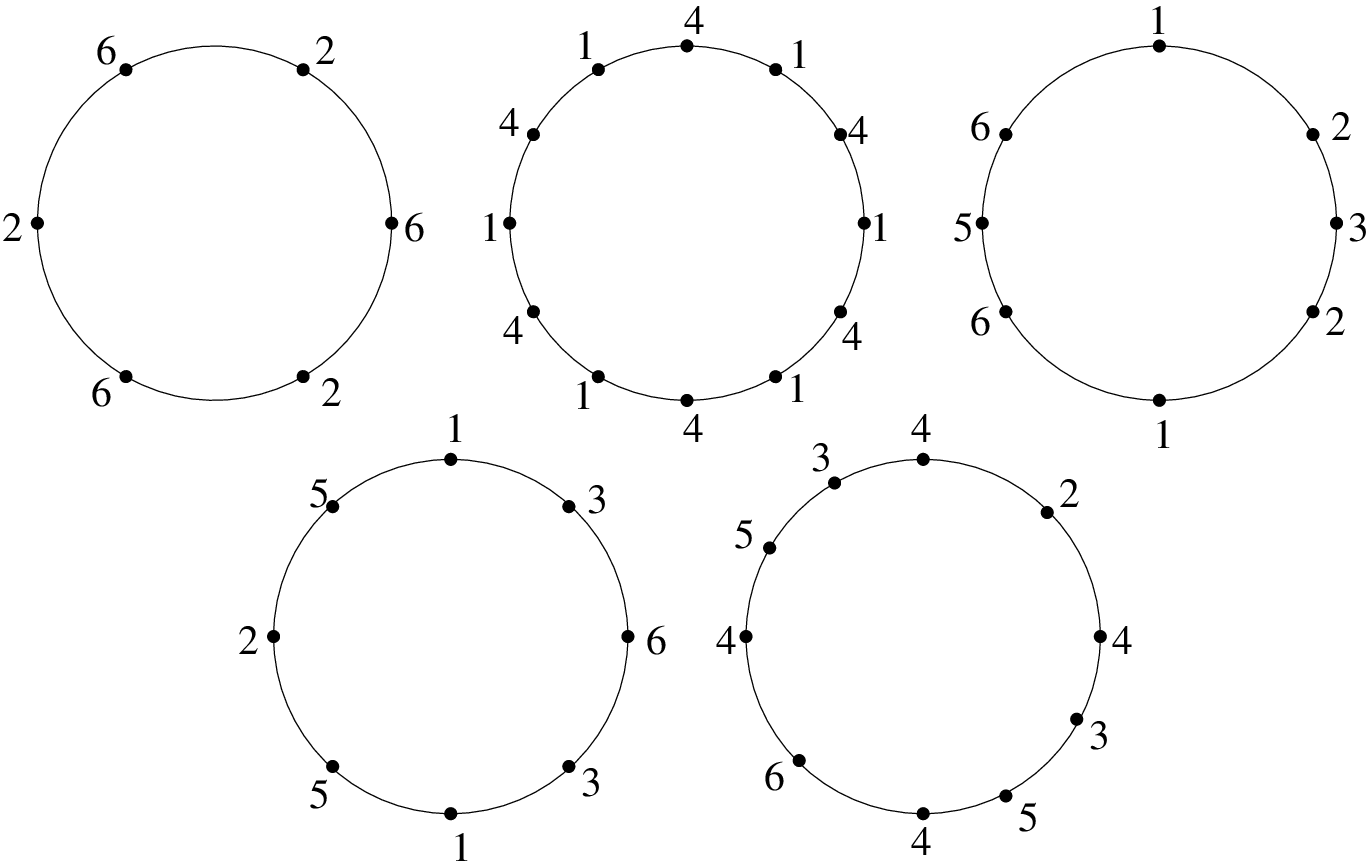}
\end{center}
We will also need some information about the geometry 
of the links of $(S^5,W_{E_6})$,
and in one case about the links of the links. 
These Coxeter complexes are of classical type, 
see section \ref{app:coxeter}. 
The links of 1-vertices have already been mentioned above. 

The link of a 2-vertex in $(S^5,W_{E_6})$ 
is a 4-dimensional Coxeter complex $\Si_2$ of type $D_5$ 
with induced labelling
\hpic{\includegraphics[scale=0.4]{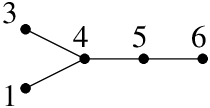}}
for its Dynkin diagram.
In view of the symmetry of the Dynkin diagram, 
the roles of the 1- and 3-vertices are equivalent. 

Any two distinct 6-vertices in $\Si_2$ 
have distance $\pihalf$ or $\pi$. 
In the first case, 
their midpoint is a 5-vertex
and they are connected by a singular segment of type $656$. 
The convex hull of a triple of 6-vertices in $\Si_2$ 
with pairwise distances $\pihalf$ 
is a right-angled equilateral triangle
centered at a 4-vertex. 
The 2-sphere containing it 
is a singular 2-sphere 
isomorphic to the $B_3$-Coxeter complex 
with Dynkin diagram
\includegraphics[scale=0.5]{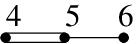}.
It is tesselated by forty eight $654$-triangles. 
\begin{center}
\includegraphics[scale=0.4]{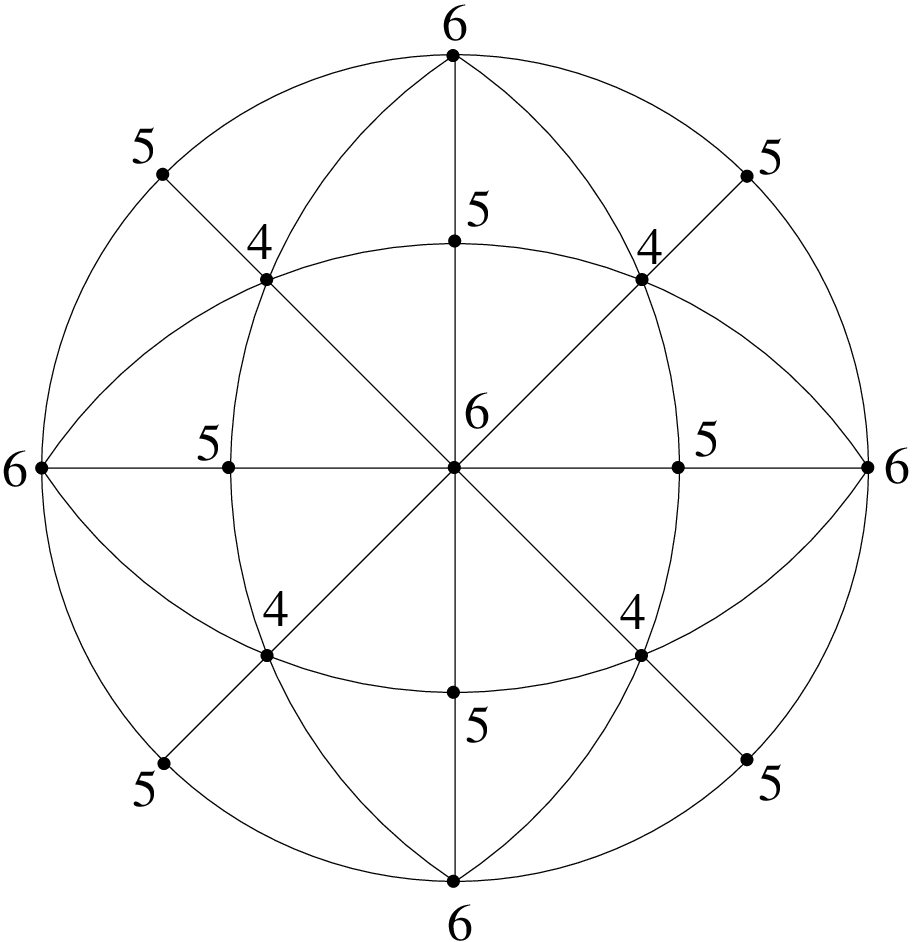}
\end{center}
The convex hull of a quadruple of 6-vertices in $\Si_2$ 
with pairwise distances $\pihalf$ 
is an equilateral tetrahedron 
with edge lengths $\pihalf$. 
Its codimension-one faces are simplicial subcomplexes 
composed of six $654$-triangles.
However, 
the tetrahedron itself is {\em not} a simplicial subcomplex;
its center is the midpoint of a $13$-edge.
Accordingly,
the geodesic 3-sphere containing the tetrahedron 
is not a subcomplex either,
and its simplicial convex hull is the entire Coxeter complex. 

The link of a 6-vertex in $(S^5,W_{E_6})$ 
is a 4-dimensional Coxeter complex $\Si_6$ of type $D_5$ 
with induced labelling
\hpic{\includegraphics[scale=0.4]{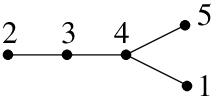}}
for its Dynkin diagram.
Its geometric properties are dual to the ones of $\Si_2$. 
For instance, 
the singular 2-sphere 
containing a $234$-triangle 
is isomorphic to the $B_3$-Coxeter complex 
with Dynkin diagram
\includegraphics[scale=0.5]{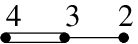}.
It is tesselated by forty eight $234$-triangles,
compare the figure for Lemma~\ref{lem:linkm3} below. 

We will use the fact 
that the possible types of 
singular segments of length $\pi$ in $\Si_6$ 
connecting antipodal 2-vertices are 
$23232$, $24342$ and $2512$.

\medskip
As already mentioned above,
the link of a 26-edge in $(S^5,W_{E_6})$ 
is a 3-dimensional Coxeter complex $\Si_{26}$ of type $D_4$ 
with induced labelling
\hpic{\includegraphics[scale=0.3]{E6link26dynk.eps}}
for its Dynkin diagram.
In view of the symmetries of the Dynkin diagram, 
the roles of the $i$-vertices for $i\neq4$ are equivalent. 

Any vertex adjacent to a 4-vertex in $\Si_{26}$ 
has distance $\piforth$ from it. 

Any two non-adjacent vertices in $\Si_{26}$
of different types $\neq4$ 
are connected by a singular segment through a vertex
of the third type $\neq4$, 
and they lie on a singular 1-sphere of type $1351351$ 
(the first and last 1 to be identified). 
For instance, 
two non-adjacent vertices of types 1 and 3 
are connected by a segment of type $153$.

Any two distinct vertices in $\Si_{26}$ of the same type $i\neq4$
have distance $\pihalf$ or $\pi$. 
In the first case, 
their midpoint is a 4-vertex
and they are connected by a singular segment of type $i4i$. 
The convex hull of a triple of $i$-vertices
with pairwise distances $\pihalf$ 
is a right-angled equilateral triangle. 
We observe that 
it is {\em not} a simplicial subcomplex of the Coxeter complex; 
its center is the midpoint of an edge perpendicular to it
and with endpoints of the two types $\neq i,4$. 
For instance, if $i=3$ then it is the midpoint of a $15$-edge. 
Accordingly, 
the geodesic 2-sphere containing such a triple of $i$-points 
is not a subcomplex, 
and its simplicial convex hull is the whole Coxeter complex. 

There is only one type of singular 2-spheres in $\Si_{26}$,
equivalently, 
the singular 2-spheres are composed of all four types 
of 2-dimensional faces. 
\begin{center}
\includegraphics[scale=0.4]{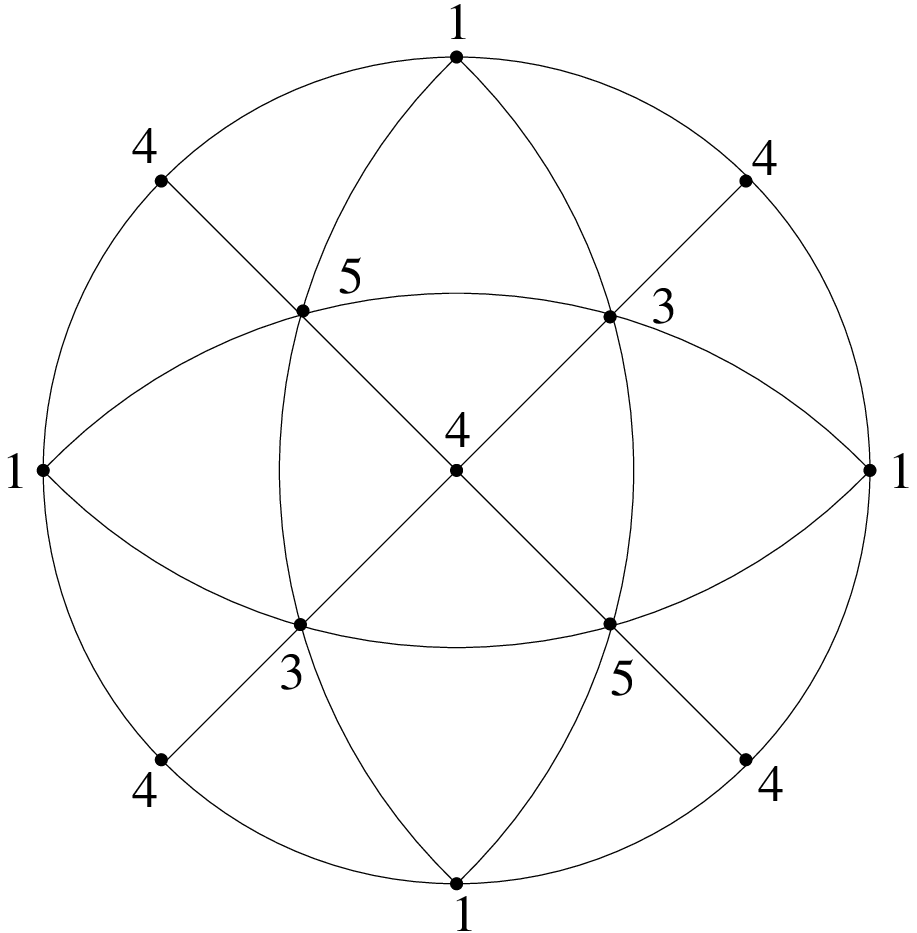}
\end{center}
\begin{rem}
\label{rem:e6max}
We have 
$Out((S^5,W_{E_6}))\cong\Z_2$. 
However,
the index two extension 
$Aut((S^5,W_{E_6}))$ of $W_{E_6}$ 
is not a reflection group,
because the nontrivial isometry of $\De$ fixes only two vertices 
and therefore is not induced by a hyperplane reflection, 
compare Remark~\ref{rem:f4max}.
\end{rem}

\subsubsection{The Coxeter complexes of classical types}
\label{app:coxeter}

In each case we consider the spherical Coxeter complex 
$(S^{n-1},W)$ as embedded in $\R^{n}$
as the unit sphere and we 
describe the Weyl group as a group of isometries of $\R^{n}$.

\medskip
{\em Type $A_n$.}
Let $n\geq 1$.
We use the labelling \includegraphics[scale=0.6]{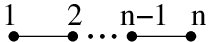} 
for the Dynkin diagram of type $A_n$.
The Weyl group $W_{A_n}$ 
is the finite group of isometries of 
$\R^n\cong \{x_0+\dots+x_n=0\}\subset \R^{n+1}$ 
generated by the reflections at the hyperplanes orthogonal 
to the {\em fundamental root vectors}
\begin{equation*}
r_i=e_{i}-e_{i-1} \text{ for } 1\leq i\leq n.
\end{equation*}
The {\em fundamental Weyl chamber} $\Delta$ 
is given by the inequalities:
\begin{equation}
\label{ineq:wcan}
x_0
\overset{\text{(1)}}{\leq} x_1
\overset{\text{(2)}}{\leq}\dots 
\overset{\text{(n)}}{\leq}x_{n}
\end{equation}
We list vectors 
representing the vertices of $\Delta$: 
\begin{center}
\begin{tabular}{R{3cm}C{1cm}R{0.1cm}R{1.6cm}R{1.6cm}R{0.7cm}R{0.7cm}R{0.7cm}R{0.9cm}R{1.1cm}}
1-vertex: & $v_1$	&	
		$($&$-n,$ & $1,$ & $1,$ & $\dots,$ & $1,$ & $1,$ & $1\;)$ \\
2-vertex: & $v_2$	&
		$($&$-(n-1),$ & $-(n-1),$ & $2,$ & $\dots,$ & $2,$ & $2,$ & $2\;)$ \\
$\vdots$ \phantom{vert} & $\vdots$	& &&&&	$\vdots$\phm &&\\
$(n-1)$-vertex: & $v_{n-1}$&
		$($&$-2,$ & $-2,$ & $-2,$ & $\dots,$ & $-2,$ & $n-1,$ & $n-1\;)$ \\
$n$-vertex: & $v_n$&	
		$($&$-1,$ & $-1,$ & $-1,$ & $\dots,$ & $-1,$ & $-1,$ & $n\;)$ \\
\end{tabular}
\end{center}
The {\em root system} consists of the vectors 
\begin{equation*}
\pm(e_i-e_j) 
\qquad\hbox{ for $0\leq i < j\leq n$.}
\end{equation*}
(The walls are the intersections of $S$ with the hyperplanes
perpendicular to a vector in the root system.) 
The Weyl group $W_{A_n}$ acts on $\R^{n+1}$ by permutations of the coordinates.

The canonical involution 
$\iota:\De_{mod}^{A_n}\to\De_{mod}^{A_n}$ 
is the nontrivial isometry. 
Accordingly, the antipodes of $i$-vertices in the Coxeter complex 
are $(n+1-i)$-vertices. 
\begin{rem}
\label{rem:anmax}
We have 
$Out((S^{n-1},W_{A_n}))\cong\Z_2$ for $n\geq2$.
However, 
the index two extension 
$Aut((S^{n-1},W_{A_n}))$ of $W_{A_n}$ 
is not a reflection group if $n\geq4$, 
because the nontrivial isometry of $\De$ moves more than two vertices 
and therefore is not induced by a hyperplane reflection. 
If $n\leq3$, 
we have 
$Aut((S^2,W_{A_3}))\cong W_{B_3}$ and $Aut((S^1,W_{A_2}))\cong W_{G_2}$. 
\end{rem}

\medskip
{\em Type $B_n$.}
Let $n\geq 2$.
We use the labelling \includegraphics[scale=0.6]{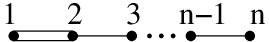} 
for the Dynkin diagram. 
The Weyl group $W_{B_n}$ 
is the finite group of isometries of 
$\R^n$ 
generated by the reflections at the hyperplanes orthogonal 
to the {\em fundamental root vectors}
\begin{equation*}
r_1=e_1,\;\;\;\; 
r_i=e_i-e_{i-1} \text{ for } 2\leq i\leq n.
\end{equation*}
The {\em fundamental Weyl chamber} $\Delta$ 
is given by the inequalities:
\begin{equation}
\label{ineq:wcbn}
0\overset{\text{(1)}}{\leq} x_1
\overset{\text{(2)}}{\leq} x_2
\overset{\text{(3)}}{\leq}\dots 
\overset{\text{(n)}}{\leq}x_n
\end{equation}
We list vectors representing the vertices of $\Delta$:
\begin{center}
\begin{tabular}{R{3cm}C{1cm}C{4cm}}
1-vertex: & $v_1$	&	$(1,1,1,\dots,1)$\\
2-vertex: & $v_2$	&	$(0,1,1,\dots,1)$\\
3-vertex: & $v_3$	&	$(0,0,1,\dots,1)$\\
$\vdots$ \phantom{vert} & $\vdots$	&	$\vdots$\\
$(n-1)$-vertex: & $v_{n-1}$ 	&	$(0,\dots,0,1,1)$\\
$n$-vertex: & $v_n$&	$(0,\dots,0,0,1)$\\
\end{tabular}
\end{center}
All half-apartments in $(S^{n-1},W_{B_n})$ are centered at a vertex. 
The vertices of types $n$ and $n-1$ 
are the vertices of {\em root type}. 
They are represented by the vectors
\begin{eqnarray*}
&\pm e_i &\qquad\hbox{ for $1\leq i\leq n$,}\\
&\pm e_i\pm e_j 
&\qquad\hbox{ for $1\leq i < j\leq n$.}
\end{eqnarray*}
The Weyl group $W_{B_n}$ acts on $\R^n$ by permutations of the coordinates
and change of signs.

The canonical involution 
$\iota:\De_{mod}^{B_n}\to\De_{mod}^{B_n}$ 
is trivial. 
Accordingly, the antipodes of $i$-vertices in the Coxeter complex 
are $i$-vertices. 
\begin{rem}
\label{rem:bnmax}
$Out((S^{n-1},W_{B_n}))$ is trivial for $n\geq3$.
For $n=2$, 
we have 
$Aut(S^1,W_{B_2})\cong W_{I_2(8)}$.
\end{rem}

\medskip
{\em Type $D_n$.}
Let $n\geq 4$. 
We use the labelling 
\hpic{\includegraphics[scale=0.5]{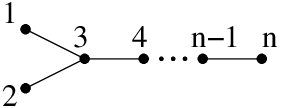}}
for the Dynkin diagram. 
The Weyl group $W_{D_n}$ 
is the finite group of isometries of 
$\R^n$ 
generated by the reflections at the hyperplanes orthogonal 
to the {\em fundamental root vectors}
\begin{equation*}
r_1=e_1+e_2,\;\;\;\; 
r_i=e_i-e_{i-1} \text{ for } 2\leq i\leq n .
\end{equation*}
The {\em fundamental Weyl chamber} $\Delta$ 
is given by the inequalities:
\begin{equation}
\label{ineq:wcdn}
-x_2\overset{\text{(1)}}{\leq} x_1
\overset{\text{(2)}}{\leq} x_2
\overset{\text{(3)}}{\leq}\dots 
\overset{\text{(n)}}{\leq}x_n
\end{equation}
We list vectors representing the vertices of $\Delta$: 
\begin{center}
\begin{tabular}{R{3cm}C{1cm}C{6cm}}
1-vertex: & $v_1$	&	$(\phm1,1,1,\dots,1)$\\
2-vertex: & $v_2$	&	$(-1,1,1,\dots,1)$\\
3-vertex: & $v_3$	&	$(\phm0,0,1,\dots,1)$\\
$\vdots$ \phantom{vert} & $\vdots$	&	$\vdots$\\
$(n-1)$-vertex: & $v_{n-1}$ 	&	$(\phm0,\dots,0,1,1)$\\
$n$-vertex: & $v_n$&	$(\phm0,\dots,0,0,1)$\\
\end{tabular}
\end{center}
All half-apartments in $(S^{n-1},W_{D_n})$ are centered at a vertex. 
The vertices of type $n-1$ 
are the vertices of {\em root type}. 
They are represented by the vectors 
\begin{equation*}
\pm e_i\pm e_j 
\qquad\hbox{ for $1\leq i < j\leq n$.}
\end{equation*}
The Weyl group $W_{D_n}$ acts on $\R^n$ by permutations of the coordinates
and change of signs in an even number of places.

The canonical involution 
$\iota:\De_{mod}^{D_n}\to\De_{mod}^{D_n}$ 
is trivial for $n$ even. 
Accordingly, in this case 
the antipodes of $i$-vertices in the Coxeter complex 
are $i$-vertices. 
For $n$ odd, the canonical involution is the nontrivial isometry.
The antipodes of 1-vertices are then 2-vertices
and for $i=3,\dots,n$ the antipodes of $i$-vertices are $i$-vertices. 

\medskip
In the models chosen here,
the  root system of $D_n$ is contained in the root system of $B_n$, 
and thus the $B_n$-Coxeter complex is a subdivision 
of the $D_n$-Coxeter complex. 
The vertices of root type ($n-1$) in $D_n$
are the vertices of a fixed type (also $n-1$) in $B_n$.
Hence the $D_n$-root system is preserved by $W_{B_n}$ and 
$W_{D_n}\subset W_{B_n}\subseteq Aut((S^{n-1},W_{D_n}))$,
cf.\ section~\ref{sec:cox}. 
If $n\geq5$, then $Out((S^{n-1},W_{D_n}))\cong\Z_2$
and it follows that 
\begin{equation}
\label{eq:autdnbn}
Aut((S^{n-1},W_{D_n}))=W_{B_n}
\qquad\hbox{ for $n\geq5$}
\end{equation}
and 
$W_{B_n}\cong W_{D_n}\rtimes\Z_2$. 
The subdivision of $\De^{D_n}$ into two $B_n$-Weyl chambers 
is obtained by cutting along the hyperplane 
perpendicular to $\pm e_1$. 

The Dynkin diagram of $\Ga_{D_4}$ has more than two symmetries,
$Aut(\Ga_{D_4})\cong Out(S^3,W_{D_4})\cong S_3$, 
and indeed the  root system of $D_4$ is contained 
in a larger root system than the one of $B_4$,
namely in the root system of $F_4$, 
Thus the $D_4$-Coxeter complex can be subdivided 
into the $F_4$-Coxeter complex. 
More precisely, 
the $D_4$-root system is preserved by $W_{F_4}$,
because the vertices of root type (3) in $D_4$
are the vertices of a fixed type (4) in $F_4$.
Accordingly,
$W_{D_4}$ is a normal subgroup of $W_{F_4}$,
$W_{D_4}\subset W_{F_4}\subseteq Aut((S^3,W_{D_4}))$.
Since also 
$W_{D_4}\subset W_{B_4}\subset W_{F_4}$,
it is clear that $[W_{F_4}:W_{D_4}]=6$.
(This can also easily be checked directly, 
for instance, 
by determining the $F_4$-walls intersecting the 
fundamental Weyl chamber $\De^{D_4}$,
respectively, containing its center.)
It follows that 
\begin{equation}
\label{eq:autd4f4}
Aut((S^3,W_{D_4}))=W_{F_4}, 
\end{equation}
and $W_{F_4}\cong W_{D_4}\rtimes S_3$. 
The subdivision of $\De^{D_4}$ into six $F_4$-Weyl chambers 
is obtained by taking the barycentric subdivision 
of the equilateral type $124$ face 
and coning off at the 3-vertex.

\subsection{Spherical buildings}
\label{sec:buil}

We refer to \cite[ch.\ 3]{qirigid}
for a treatment of spherical buildings 
from the perspective of comparison geometry. 

\subsubsection{Some basic definitions}

We recall the geometric definition of spherical buildings 
as given in \cite[ch.\ 3.2]{qirigid}. 
A {\em spherical building} 
modelled on a Coxeter complex $(S,W)$ 
is a CAT(1) space $B$ 
equipped with an {\em atlas} 
consisting of isometric embeddings $\iota:S\embed B$, 
the {\em charts}, 
satisfying certain properties. 
The images of the charts are called the {\em apartments}.
Any two points must lie in an apartment. 
The atlas must be closed under precomposition 
with isometries in $W$,  
and the charts must be compatible in the sense that 
the coordinate changes are restrictions of isometries in $W$. 

The underlying set may be empty, 
in which case the building is called a spherical {\em ruin}. 

One defines walls, roots, singular spheres, 
faces, chambers, panels, 
regular and singular points 
as the images of the corresponding objects in the model Coxeter complex. 
In particular, 
a spherical building carries a natural structure 
as a piecewise spherical {\em polyhedral} complex, 
in fact, as a {\em simplicial} complex 
if $W$ has no fixed points on $S$. 
The building $B$ is called {\em thick}, 
if every panel is adjacent to at least three chambers. 

Two points in $B$ are called {\em antipodal} 
if they have maximal distance $\pi$.

Let $\si\subset B$ be a face of codimension $\geq1$. 
Then for an interior point $x\in\si$, 
the link $\Si_xB$ splits as the spherical join 
$\Si_xB\cong\Si_x\si\circ\nu_x\si$
of the unit sphere $\Si_x\si$ 
and the unit normal space $\nu_x\si$ of $\si$ in $B$. 
One can consistently identify with each other 
the unit normal spaces $\nu_x\si$ for all interior points $x\in\si$. 
This identification can be described as follows:  
For interior points $x_1,x_2\in\si$, 
let $c_i:[0,\eps)\to B$ 
be unit speed geodesic segments 
emanating from $x_i$ orthogonal to $\si$. 
Then the directions $\dot c_i(0)\in\nu_{x_i}\si$ 
are identified if and only if for small $t>0$ 
the convex hulls $CH(\si\cup\{c_i(t)\})$ 
locally coincide near $x_1$ and $x_2$.  
We call the resulting identification space 
the {\em link $\Si_{\si}B$ of the face $\si$}. 
It inherits a natural structure 
as a spherical building modelled on the Coxeter complex 
$\bigl(\Si_{\iota^{-1}(\si)}S,Stab_W(\iota^{-1}(\si))\bigr)$
where $\iota:S\embed B$ is a chart with $\si\subset\iota(S)$. 
For faces $\si\subset\tau$ 
there is a canonical identification 
$\Si_{\tau}B\cong\Si_{\Si_{\si}\tau}(\Si_{\si}B)$.

There is a natural ``accordion'' {\em anisotropy map}
$\theta_B:B\to\De_{mod}$ 
folding the building onto the model Weyl chamber. 
It is determined by the property 
that for any chart $\iota:S\embed B$ holds 
$\theta_B\circ\iota=\theta_S$. 
The anisotropy map is 1-Lipschitz and restricts to isometries on chambers. 
Furthermore, 
for any apartment $a\subset B$ and any chamber $\si\subset a$ 
there is a natural 1-Lipschitz retraction 
$\rho_{a,\si}:B\to a$, 
$\rho_{a,\si}|_a=id_a$, 
which can be described as follows: 
For a regular point $x\in\si$ 
and any point $y\in B$ not antipodal to $x$, 
the segment connecting $x$ to $\rho_{a,\si}(y)$
coincides with the segment $xy$ near $x$ and has the same length. 
The retraction is type preserving, 
$\theta_B\circ\rho_{a,\si}=\theta_B$, 
and it restricts to an isometry 
on every apartment containing $\si$. 

The following result in the spirit of 
\cite[Prop.\ 3.5.1]{qirigid} 
allows to recover the building structure from the anisotropy map. 

\begin{prop}
[Recognizing a building structure]
\label{prop:recog}
Suppose that $X$ is a CAT(1) space 
which is equipped with a structure as a 
piecewise spherical polyhedral complex 
of dimension equal to $\dim(S)$. 
Let $\theta_X:X\to\De_{mod}$ 
be a 1-Lipschitz map 
which restricts on every top-dimensional face 
to an isometry onto $\De_{mod}$. 
Suppose furthermore that any two points in $X$ 
lie in an isometrically embedded copy of $S$. 
Then $X$ carries a natural structure as a spherical building 
modelled on the Coxeter complex $(S,W)$
and with anisotropy map $\theta_X$. 
\end{prop}
\proof 
We call a top-dimensional face of $X$ a {\em chamber}
and an isometrically embedded copy of $S$ an {\em apartment}. 
Due to our assumption, 
the apartments are tesselated by chambers. 
For any two adjacent chambers $\si_1$ and $\si_2$ 
in an apartment $a$ 
there is the isometry 
$(\theta_X|_{\si_2})^{-1}\circ(\theta_X|_{\si_1}):\si_1\to\si_2$. 
It must coincide with the reflection at the common 
codimension-one face $\si_1\cap\si_2$. 
Hence there is an isometric identification $\iota:S\to a$ 
satisfying $\theta_X\circ\iota=\theta_S$ 
which is unique up to precomposition with Weyl isometries. 
The compatibility of all these {\em charts} for all apartments 
is automatic 
and they form an atlas for a spherical building structure 
modelled on $(S,W)$ 
with anisotropy map $\theta_X$. 
\qed

An isometry $\al:B\to B$ is called an {\em automorphism}
of the spherical building $B$ 
if  it preserves the polyhedral structure. 
If $B$ is a {\em thick} building,
then all its isometries are automorphisms. 
An automorphism is called {\em inner} 
if it is type preserving, 
$\theta_B\circ\al=\theta_B$. 
We denote the automorphism group of $B$ by $Aut(B)$, 
and the subgroup of inner automorphisms by $Inn(B)$. 
Then $Inn(B)$ is a finite index normal subgroup of $Aut(B)$ 
and $Out(B):=Aut(B)/Inn(B)$ embeds as a subgroup of 
$Isom(\De_{mod})\cong Out((S,W))$.

\subsubsection{Convex subcomplexes and subbuildings}
\label{sec:subcomp}

Let $B$ be a spherical building. 
By a {\em convex subcomplex} of $B$ 
we mean a closed convex subset 
which is also a subcomplex
with respect to the natural polyhedral structure on $B$.
The {\em simplicial convex hull} of a subset $A\subseteq B$ 
is the smallest convex subcomplex of $B$ containing $A$. 

We call a convex subcomplex $K\subseteq B$ 
a {\em subbuilding} 
if any two of its points are contained in a singular sphere $s\subseteq K$ 
with $\dim(s)=\dim(K)$. 
The next result tells that 
a subbuilding inherits a natural structure as a spherical building. 
To describe the associated Coxeter complex, 
let $a\subset B$ be an apartment containing 
a singular sphere $s\subset K$ with $\dim(s)=\dim(K)$  
and let 
$\iota:S\buildrel\cong\over\to a\subset B$ 
be a chart. 
As explained in section~\ref{sec:cox}, 
the singular sphere $\iota^{-1}(s)\subseteq S$ 
inherits from $S$ a natural structure as a Coxeter complex
with a possibly coarser polyhedral structure. 
Its {\em $W$-type}, that is, 
its equivalence class modulo the action of the Weyl group 
does not depend on the choice of $s$. 

\begin{prop}
[Building structure on subbuildings]
\label{prop:subbuil}
The subbuilding $K\subseteq B$ 
carries a natural structure as a spherical building 
modelled on the Coxeter complex
$(\iota^{-1}(s),W_{\iota^{-1}(s)})$. 
\end{prop}
\proof
We fix a singular sphere $s\subset K$ with $\dim(s)=\dim(K)$. 
Let $a$ be an apartment containing $s$ 
and let $\si\subset a$ be a chamber 
such that $\si\cap s$ is a top-dimensional face of $s$. 
Then the retraction 
$\rho_{a,\si}:B\to a$
restricts to a retraction 
$\rho_{s,\si\cap s}:K\to s$ of $K$.

We note that $\rho_{s,\si\cap s}$ 
restricts to an isometry 
on every singular sphere $s'\subset K$ containing $\si\cap s$. 
Using a chart 
$\iota:S\buildrel\cong\over\to a\subset B$, 
we can pull back the intrinsic polyhedral structure 
(as a Coxeter complex) on $\iota^{-1}(s)$ to $s'$ 
via $\iota^{-1}\circ\rho_{s,\si\cap s}|_{s'}$. 
We will refer to the pulled back structure 
as the {\em intrinsic} polyhedral structure on $s'$. 
The main point to verify is 
that the intrinsic polyhedral structures 
on all such singular spheres $s'$ match 
and yield a polyhedral structure on $K$. 

At this point, we have on $K$ only the polyhedral structure 
which it inherits from $S$. 
We say that $K$ {\em branches}
along a codimension one face $\phi$, 
if $K$ has at least three 
top-dimensional faces $\tau_1,\tau_2$ and $\tau_3$ 
adjacent to $\phi$, 
i.e.\ with $\phi$ as a common codimension one face. 
It then follows that the $\tau_i$ 
(and all top-dimensional faces of $K$ adjacent to $\phi$) 
must have the same $\theta_B$-type. 
Indeed, 
the unions $\tau_i\cup\tau_j$ for $i\neq j$ are convex, 
because $K$ is convex and the $\tau_i$ are top-dimensional in $K$. 
Let $\iota_{ij}$ be charts 
whose images contain $\tau_i\cup\tau_j$. 
We may choose, say, $\iota_{12}$ and $\iota_{13}$
so that $\iota_{12}^{-1}(\tau_1)=\iota_{13}^{-1}(\tau_1)$. 
Then necessarily 
$\iota_{12}^{-1}(\tau_2)=\iota_{13}^{-1}(\tau_3)$ 
and hence 
$\theta_B(\tau_2)=\theta_B(\tau_3)$. 
Similarly, 
$\theta_B(\tau_1)=\theta_B(\tau_2)$. 

Let $s'\subset K$ be a singular sphere 
containing $\phi$ and $\si\cap s$. 
Since the two top-dimensional faces in $s'$ adjacent to $\phi$ 
have the same $\theta_B$-type, 
the codimension one singular subsphere $t'\subset s'$ containing $\phi$ 
is an {\em $s'$-wall}
in the sense that 
$(\iota^{-1}\circ\rho_{s,\si\cap s})(t')$ 
is an $\iota^{-1}(s)$-wall 
as defined in section~\ref{sec:cox}.

The intersection 
$s'\cap s''$ 
of any two singular spheres 
$s',s''\subset K$ 
containing $\si\cap s$ 
is obviously top-dimensional in $K$. 
Our discussion implies that it is a subcomplex 
with respect to the {\em intrinsic} polyhedral structures 
of these singular spheres, 
because its boundary consists of codimension one faces of $K$
along which $K$ branches 
and thus is contained in a union of $s'$-walls ($s''$-walls).
It follows that the intersection of 
any two top-dimensional faces 
$\tau'\subset s'$ and $\tau''\subset s''$ 
with respect to the intrinsic polyhedral structures 
is either empty or a face of $\tau'$ (and $\tau''$). 
This means that the intrinsic polyhedral structures 
on the singular spheres in $K$ containing $\si\cap s$ 
match and form together a polyhedral structure on $K$. 

To conclude the argument,
we observe that the 1-Lipschitz map
\begin{equation*}
\theta_{\iota^{-1}(s)}\circ\iota^{-1}\circ\rho_{s,\si\cap s}:K
\to\iota^{-1}(s)/W_{\iota^{-1}(s)}
=\De_{mod}^{(\iota^{-1}(s),W_{\iota^{-1}(s)})}
\end{equation*}
restricts 
on top-dimensional faces 
(for the new polyhedral structure just defined)
to surjective isometries,
because for singular spheres $s'\subset K$ 
containing $\si\cap s$
the map 
$\iota^{-1}(s)\circ\rho_{s,\si\cap s}|_{s'}:s'\to\iota^{-1}(s)$ 
is an isometry preserving the polyhedral structure. 
The assertion follows now from  
Proposition~\ref{prop:recog}. 
\qed

\medskip
We discuss now some conditions 
implying that a convex subcomplex is a subbuilding. 
To begin with,
the existence of a top-dimensional subsphere is sufficient. 
\begin{prop}
\label{prop:contsphsubbuil}
If a convex subcomplex $K\subseteq B$ 
contains a singular sphere $s\subseteq K$ with $\dim(s)=\dim(K)$, 
then it is a subbuilding. 
\end{prop}
\proof
For top-dimensional subcomplexes
this is \cite[Prop.\ 3.10.3]{qirigid}. 
The proof in the general case is similar. 
One observes first that 
every point $x\in K$ has an antipode $\hat x$ in $s$. 
More precisely, 
for an arbitrary point $y\in s$ 
the antipode $\hat x$ can be chosen 
so that $y$ lies on a geodesic segment $x\hat x$ of length $\pi$. 
To see this,
let $\si\subset s$ be a face containing $y$ with $\dim(\si)=\dim(s)$
and let $y_1$ be an interior point of $\si$.
Then the geodesic segment $xy_1$
is contained in $s$ near $y_1$
and therefore can be extended beyond $y_1$ inside $s$
to a geodesic segment $xy_1\hat x$ of length $\pi$. 
The convex hull $CH(\{x,\hat x\}\cup\si)$ is a bigon 
and contains a geodesic segment $xy\hat x$ of length $\pi$.

The convex hull of $x$ and a small disk in $s$ around $\hat x$ 
is a singular sphere $s'$ with $y\in s'\subset K$. 
We see that any point and also any pair of points in $K$ 
lies in a singular sphere of dimension $\dim(K)$ contained in $K$. 
\qed

If every point of $K$ has an antipode {\em in} $K$, 
then clearly plenty of spheres $s\subset K$ 
as in Proposition~\ref{prop:contsphsubbuil} exist. 
Just take the simplicial convex hull of a pair 
of maximally regular antipodes in $K$. 
(A point in $K$ is said to be 
{\em maximally regular} 
if it is an interior point of a face $\si\subset K$ 
which is top-dimensional in $K$, 
$\dim(\si)=\dim(K)$.)
The next result says that it is enough 
to assume the existence of antipodes only for vertices. 

\begin{prop}[{\cite[Thm.\ 2.2]{Serre}}]
\label{lem:allvertantip}
If every vertex of a convex subcomplex $K\subseteq B$ 
has an antipode in $K$, 
then $K$ is a subbuilding. 
\end{prop}
\proof
Let $\si\subset K$ be a simplex with vertices $p_0,\dots,p_k$. 
In view of Proposition~\ref{prop:contsphsubbuil},
it suffices to show that there exists a singular $k$-sphere $s_k$
with $\si\subset s_k\subset K$.

We proceed by induction over $k$. 
By assumption, the assertion holds for $k=0$. 
To do the induction step,
we consider a singular $(k-1)$-sphere $s_{k-1}$ with 
$p_0,\dots,p_{k-1}\in s_{k-1}\subset K$.
The convex hull of $s_{k-1}$ and $p_k$ is a $k$-hemisphere 
$h_k\subset K$ 
with boundary $\partial h_k=s_{k-1}$. 
By assumption, $p_k$ has an antipode $\hat p_k$ in $K$. 
The convex hull of $h_k$ and $\hat p_k$ contains 
a singular $k$-sphere $s_k$ 
with $\si\subset s_k\subset K$. 
\qed

\medskip
The following simple observation will be useful 
when we search for antipodes in convex subcomplexes.
\begin{lem}
\label{lem:antipodes}
Let $x_1x_2\subset K$ be a segment. 
Suppose that $z$ is an interior point of $x_1x_2$ 
which has an antipode $\hat z\in K$. 
Then the $x_i$ also have antipodes in $K$.
\end{lem}
\proof
Let $\gamma_i\subset K$ for 
$i=1,2$ be the geodesic connecting $z$ and $\hat z$
with initial direction $\ora{zx_{3-i}}$ at $z$. 
Then $x_iz \cup \gamma_i$ is a geodesic of length $>\pi$, 
and $\gamma_i$ contains an antipode of $x_i$.
\qed

\medskip
We will call a point $x$ in a convex subcomplex $K\subset B$
an {\em interior} point of $K$ 
if $\Si_xK$ is a subbuilding of $\Si_xB$,
and a {\em boundary} point otherwise.

\subsubsection{Circumcenters}

Let $B$ be a spherical building 
and let $K\subset B$ be a non-empty convex subcomplex. 

We recall that 
if $\rad(K)<\pihalf$ 
then $K$ has a unique circumcenter 
which must be contained in $K$, 
cf.\ section~\ref{sec:rad}. 

If $\rad(K)=\pihalf$, 
then the arguments in \cite[ch.\ 3]{BalserLytchak_centers} 
for general CAT(1) spaces of finite dimension 
yield that $Cent(K)\neq\emptyset$. 
Moreover, 
if $K\cap Cent(K)\neq\emptyset$, 
then $K\cap Cent(K)$ has a unique circumcenter. 
In our special situation, 
the proofs simplify 
and we include them for the sake of completeness. 
\begin{lem}
\label{lem:centex}
If $B$ is a spherical building 
and if $K\subset B$ is a convex subcomplex 
with $\rad(K)=\pihalf$, 
then $Cent(K)\neq\emptyset$. 
\end{lem}
\proof
Consider a sequence of points $x_i\in B$ 
with $\rad(K,x_i)\searrow\pihalf$. 
We may assume that $\theta_B(x_i)\to t\in\De_{mod}$
by compactness. 
Let $\si_i$ denote the face of $B$ 
containing $x_i$ as an interior point. 
For sufficiently large $i$, 
$\si_i$ contains a point $x'_i$ with $\theta_B(x'_i)=t$. 
(Namely, when $d(\theta_B(x_i),t)<\de$
where $\de$ is the distance between $t$ 
and the union of those faces of $\De_{mod}$ which do not contain $t$ 
in their closure.) 
We have $d(x'_i,x_i)\to0$ and hence also $\rad(K,x'_i)\to\pihalf$. 
Since the $x'_i$ have fixed $\theta_B$-type $t$
and since $K$ is a subcomplex, 
it follows that $\rad(K,x'_i)=\pihalf$ for large $i$. 
This is due to the fact that 
the radius of a face of $B$ 
with respect to a point of fixed type $t$ 
can take only finitely many values 
(depending on $t$ and the Coxeter complex). 
\qed

\medskip
The following observations apply to the situation 
when the closed convex subset $K\cap Cent(K)$ is non-empty. 
Clearly, it has diameter $\leq\pihalf$. 
We begin with a consequence of \cite[Prop.\ 1.2]{BalserLytchak_centers}:
\begin{lem}
\label{lem:diampihalf}
Suppose that $A\subset B$ is a non-empty subset 
with $\diam(A)\leq\pihalf$ 
and such that $\theta_B(A)\subset\De_{mod}$ is finite. 
Then $\rad(A)<\frac{\pi}{2}$. 
\end{lem}
\proof
By the finiteness assumption on types,
the distances between points in $A$ take only finitely many values. 
Therefore, 
if for some point $x\in A$ holds $d(x,y)<\pihalf$ for all $y\in A$, 
then $\rad(A,x)<\pihalf$ and we are done. 

Otherwise, 
we pick some point $x\in A$
and consider the set $A'\subset\Si_xB$ of directions $\ora{xy}$
to the points $y\in A$ with $d(x,y)=\pihalf$. 
The directions in $A'$ have only finitely many 
$\theta_{\Si_xB}$-types and, 
by triangle comparison,
$\diam(A')\leq\pihalf$.
Thus $A'$ satisfies the same assumptions as $A$. 
Moreover, 
$\rad(A')<\pihalf$ implies $\rad(A)<\pihalf$. 
Indeed, 
let $\ga:[0,\eps)\to B$ be a geodesic segment 
with initial point $\ga(0)=x$ and initial direction $\dot\ga(0)$
satisfying $\rad(A',\dot\ga(0))<\pihalf$.
Then $\rad(A,\ga(t))<\pihalf$ for small $t>0$. 
Here we use again the finiteness of $\theta_B(A)$;
it yields a constant $\eps>0$ with the property that 
$d(x,y)=\pihalf$ or $d(x,y)\leq\pihalf-\eps$ 
for all $y\in A$. 

We can therefore proceed by induction on the dimension of $B$. 
The assertion holds trivially for buildings of dimension zero. 
\qed

\begin{cor}
\label{cor:diampihalffixpt}
Let $C\subset B$ be a non-empty closed convex subset 
with diameter $\leq\pihalf$. 
Then the action $Stab_{Aut(B)}(C)\acts C$ has a fixed point. 
\end{cor}
\proof
Pick any point in $C$.
By Lemma~\ref{lem:diampihalf}, 
its $Stab_{Aut(B)}(C)$-orbit has circumradius $<\pihalf$ 
and therefore a unique circumcenter which is contained in $C$.
It is fixed by $Stab_{Aut(B)}(C)$. 
\qed

\begin{cor}
\label{cor:fixptoncompl}
If $K$ is as in Lemma~\ref{lem:centex} 
and if $K\cap Cent(K)\neq\emptyset$, 
then $Stab_{Aut(B)}(K)\acts K$ has a fixed point. 
\end{cor}
\proof
By Corollary~\ref{cor:diampihalffixpt},
$Stab_{Aut(B)}(K)$ fixes a point in $K\cap Cent(K)$.
\qed

\begin{rem}
\label{rem:fixedpt}
$Cent(K)$ is the intersection of 
the closed $\pihalf$-balls centered at the vertices of $K$. 
Hence $Cent(K)$ is a subcomplex of $B$ 
with respect to a refinement of the polyhedral structure of $B$ 
which corresponds to a refinement of the polyhedral structure 
of the Coxeter complex $(S,W)$. 
This refinement can be described as follows: 
Consider the boundaries of the hemispheres 
in the Coxeter complex centered at its vertices.
In general, not all of these codimension one great spheres 
are walls. 
However there are only finitely many of them
and they yield a refinement of the polyhedral structure of $(S,W)$ 
which projects to a subdivision of $\De_{mod}$. 
If $K\cap Cent(K)\neq\emptyset$, 
we may therefore apply Lemma~\ref{lem:diampihalf} 
to the set of vertices of $K\cap Cent(K)$ 
with respect to the refined polyhedral structure 
and conclude that $\rad(K\cap Cent(K))<\pihalf$. 
It follows that $K\cap Cent(K)$ contains a unique circumcenter 
which must be fixed by $Isom(K)\supseteq Stab_{Aut(B)}(K)$.
\end{rem}

Combining these results with Proposition~\ref{prop:contsphsubbuil}, 
we obtain: 

\begin{cor}
\label{cor:cod1sph}
If a convex subcomplex $K\subseteq B$ 
contains a singular sphere $s\subseteq K$ with $\dim(s)=\dim(K)-1$, 
then $K$ is a subbuilding 
or it is contained in a closed $\pihalf$-ball 
centered in $K$ 
and the action $Stab_{Aut(B)}(K)\acts K$ has a fixed point. 
\end{cor}
\proof
$K$ contains a (polyhedral) hemisphere $h\subseteq K$ with boundary $s$. 
One can obtain $h$ as the convex hull of $s$ 
and a top-dimensional face $\si$ of $K$ 
such that $\si\cap s$ is a top-dimensional face of $s$. 
Let $z$ be the center of $h$. 

If $K\subseteq\ol B_{\pihalf}(z)$, 
then $\rad(K)=\pihalf$ 
and $z\in K\cap Cent(K)$.
Corollary~\ref{cor:fixptoncompl} yields 
that $Stab_{Aut(B)}(K)\acts K$ has a fixed point. 

If $K\not\subseteq\ol B_{\pihalf}(z)$, 
let $x\in K$ be a point with $d(x,z)>\pihalf$. 
There exists an antipode $\hat x\in h-s$ of $x$. 
(Just pick a maximally regular point $y\in h$ close to $z$ 
and extend the geodesic segment $xy$ beyond $y$ inside $h$
up to length $\pi$.) 
It follows that 
the convex hull of $x$ and $h$ contains a singular sphere $s'$ 
with $\dim(s')=\dim(K)$. 
Proposition~\ref{prop:contsphsubbuil} 
then implies that $K$ is a subbuilding. 
\qed

\medskip
The proof shows that 
if $K$ is not a subbuilding
and if $h\subset K$ is a hemisphere with $\dim(h)=\dim(K)$ and center $z$,
then $K\subseteq\ol B_{\pihalf}(z)$.

\section{On the Center Conjecture}
\subsection{General properties of potential counterexamples}
\label{sec:counter}

Let $B$ be a spherical building 
and let $K\subseteq B$ be a convex subcomplex. 
We call $K$ a {\em counterexample}
to the Center Conjecture~\ref{conj:tcc},
if $K$ is not a subbuilding 
and if the action $Stab_{Aut(B)}(K)\acts K$ 
has no fixed point. 

It is easy to see 
that a one-dimensional convex subcomplex $G\subset B$
is either a subbuilding or a metric tree 
with intrinsic circumradius $\leq\pihalf$. 
In the latter case, 
$G\cap Cent(G)$ consists of precisely one point
which then must be fixed by the action $Isom(G)\acts G$. 

Hence a counterexample $K$ must have dimension $\geq2$. 
According to \cite{BalserLytchak_centers}, 
it must even have dimension $\geq3$. 
However, we will not use this fact 
in order to keep our arguments self-contained.

Our considerations in chapter~\ref{sec:prelim} imply 
that $K$ cannot contain a singular sphere of codimension one (in $K$), 
cf.\ Corollary~\ref{cor:cod1sph}. 
Furthermore, 
$K$ can neither contain 
a $Stab_{Aut(B)}(K)$-invariant subset with circumradius $<\pihalf$ 
nor one with diameter $\leq\pihalf$, 
cf.\ Corollary~\ref{cor:diampihalffixpt}.

\subsection{The $F_4$-case}
\label{sec:f4}

We now prove our first main result.
\begin{thm}
\label{thm:tccf4}
The 
Center Conjecture~\ref{conj:tcc} 
holds for spherical buildings of type $F_4$. 
\end{thm}
\proof
We will use the information 
on the geometry of the $F_4$-Coxeter complex 
collected in section~\ref{sec:f4geom}.

Let $B$ be a spherical building of type $F_4$ 
and let $K\subseteq B$ be a convex subcomplex 
which is a counterexample in the sense of 
section~\ref{sec:counter}. 
Then $K$ must have dimension 2 or 3. 

We start by checking that also the action 
of the potentially smaller group of {\em inner} automorphisms of $B$ 
preserving $K$ has no fixed point on $K$. 
\begin{lem}
\label{lem:f4innnofix}
The action $Stab_{Inn(B)}(K)\acts K$ has no fixed point. 
\end{lem}
\proof
We assume that there exists an automorphism 
$\al\in Aut(B)-Inn(B)$ preserving $K$,
because otherwise there is nothing to prove.
In view of the natural embedding 
$Out(B)\embed Isom(\De_{mod}^{F_4})\cong\Z_2$, 
$\al$ induces on $\De_{mod}^{F_4}$
the nontrivial isometric involution. 
Hence it switches the vertex types
$1\leftrightarrow4$ and $2\leftrightarrow3$. 
Moreover,
$Stab_{Aut(B)}(K)$ is generated by $\al$ 
and the index two normal subgroup $Stab_{Inn(B)}(K)$. 

The fixed point set of any inner automorphism of $B$ 
is a convex subcomplex, 
therefore also the $Stab_{Aut(B)}(K)$-invariant subset 
$F:=K\cap Fix(Stab_{Inn(B)}(K))$. 
Note that $\al$ acts on $F$ as an isometric involution 
without fixed point (since $K$ is a counterexample). 
Hence $\al$ must map any point $x\in F$ to an antipode, 
because otherwise $x$ and $\al(x)$ would have a unique midpoint 
which would be a fixed point of $\al$ in $F$.
On the other hand, 
for any vertex $v\in F$, 
$v$ and $\al(v)$ have different $\theta_B$-types
and therefore cannot be antipodal. 
This shows that $F=\emptyset$, as claimed.
\qed

\medskip\noindent
{\em Proof of Theorem~\ref{thm:tccf4} continued.}
Due to the symmetry of the Dynkin diagram, 
the roles of $i$- and $(5-i)$-vertices are equivalent
(dual). 
Our strategy will be 
to investigate the pattern of 1- and 4-vertices in $K$, 
because vertices of these types are better separated from each other 
than 2- and 3-vertices.
We recall that the possible mutual distances between 1-vertices 
(4-vertices) in $B$ are 
$0,\pithird,\pihalf,\2pithird$ and $\pi$. 

{\em Case 1:
All vertices of types 1 and 4 in $K$ have antipodes in $K$.}
Let $q\in K$ be a 2-vertex. 
Since $dim(K)\geq2$, there
is a vertex $p\in K$ of type 1 or 4 adjacent to $q$.
Let $\hat p\in K$ be an antipode of $p$.
The geodesic segment $pq\hat p$ of length $\pi$
contains an interior 1-vertex $p'$ (between $q$ and $\hat p$).
Since also $p'$ has an antipode in $K$,
it follows from Lemma~\ref{lem:antipodes} that $q$ has an antipode in $K$.
Hence all 2-vertices in $K$ have antipodes in $K$, 
and analogously all 3-vertices in $K$ do. 
This contradicts Proposition~\ref{lem:allvertantip}.

{\em Case 2:
$K$ contains vertices of type 1 or 4 without antipodes in $K$.}
We may assume without loss of generality
that $K$ contains 1-vertices without antipodes in $K$.

If $P$ is any property defined for $i$-vertices in $K$
and invariant under 
$Stab_{Inn(B)}(K)$, 
we call vertices with this property {\em $iP$-vertices}.
If $1P$-vertices exist, 
then for any $1P$-vertex there is another $1P$-vertex 
at distance $\2pithird$ or $\pi$. 
Namely, 
let $1P'$ be property $1P$ with the additional requirement 
that all other $1P$-vertices have distance $\leq\pihalf$. 
The set of $1P'$-vertices has diameter $\leq\pihalf$ 
and must therefore be empty 
in view of Corollary~\ref{cor:diampihalffixpt} 
and Lemma~\ref{lem:f4innnofix}. 

Let $A$ be the property of {\em not having antipodes in $K$}. 
If $P$ is a property implying $A$, 
$P\Rightarrow A$,  
then for any $1P$-vertex 
there exists another $1P$-vertex at distance {\em exactly} $\2pithird$.

Let $I$ be the property of {\em being an interior point of $K$},
compare the definition at the end of section~\ref{sec:subcomp}. 
According to Proposition~\ref{prop:contsphsubbuil}, 
a point $x\in K$ is an interior point of $K$, 
if and only if $\Si_xK$ contains a (singular) sphere of top dimension 
$\dim(K)-1$. 
Clearly $I\Rightarrow A$, 
because $K$ contains no top-dimensional (in $K$) singular spheres. 

Let $x_1$ and $x_2$ be $1I$-vertices 
with distance $\2pithird$. 
The segment $x_1x_2$ is of type $12121$
and can be extended inside $K$ beyond both endpoints
to a segment $y_1y_2$ of length $\pi$ and type $2121212$. 
Note that the links of 2-vertices in $B$ 
have the type $A_1\circ A_2$ Dynkin diagram 
\includegraphics[scale=0.6]{F4link2dynk.eps}.
Since $x_i$ is an interior point of $K$,
the link $\Si_{y_i}K$ contains a top-dimensional (in $\Si_{y_i}K$) hemisphere 
centered at $\ora{y_ix_i}$. 
Consequently, $K$ contains a top-dimensional (in $K$) hemisphere 
(with $y_1$ and $y_2$ in its boundary). 
This is impossible 
and we see that $K$ cannot contain $1I$-vertices.
Dually, $K$ cannot contain $4I$-vertices. 

Below,
we will consider yet another property.
Note that a singular 1-sphere of type $\dots4343\dots$ 
in the $B_3$-Coxeter complex with Dynkin diagram 
\includegraphics[scale=0.6]{F4link1dynk.eps}
divides it into two hemispheres centered at 4-vertices. 
We say that a 1-vertex $x\in K$ has property $H$,
if {\em $\Si_xK$ contains a 2-dimensional hemisphere 
centered at a 4-vertex}. 
Let $u\in K$ denote the 4-vertex adjacent to $x$ 
such that $\ora{xu}$ is the center of this hemisphere. 
We note that $\Si_xK\subseteq\ol B_{\pihalf}(\ora{xu})$
by (the proof of) Corollary~\ref{cor:cod1sph}, 
because $K$ contains no $1I$-vertices. 
One gives a dual definition of property $H$ for 4-vertices of $K$. 
Clearly, $H\Rightarrow A$ 
for 1- and 4-vertices 
because $K$ contains no 3-dimensional hemisphere.

\parpic{\includegraphics[scale=0.7]{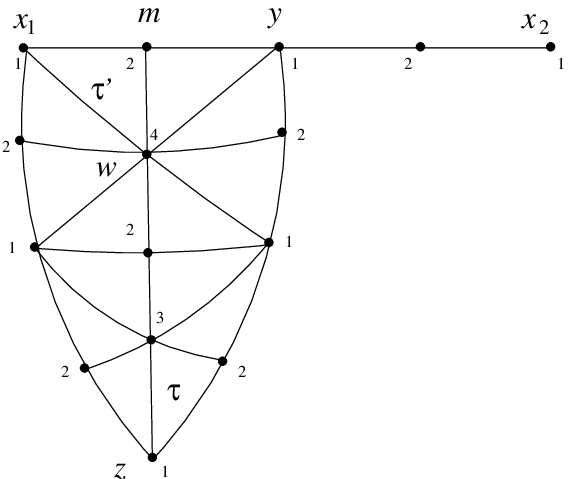}}
Let us return to an arbitrary property $P$ implying $A$. 
Let $x_1$ and $x_2$ be a pair of $1P$-vertices 
with distance $\2pithird$. 
The midpoint $y$ of $x_1x_2$ is a $1A$-vertex 
(cf. Lemma~\ref{lem:antipodes}). 
There exists another 1-vertex $z\in K$ with $d(y,z)=\2pithird$.
Both angles $\angle_y(z,x_i)$ are $<\pi$, 
because $d(z,x_i)<\pi$. 
At least one of them, say $\angle_y(z,x_1)$, must be $\geq\pihalf$
and the arc in $\Si_yK$ 
connecting $\ora{yz}$ with $\ora{yx_1}$
is then of type $242$.
We observe that the convex hull $CH(x_1,y,z)$
is an isosceles spherical triangle 
(meaning the two-dimensional object). 
Indeed, 
let $m$ denote the midpoint of the segment $yx_1$. 
The convex hull of $ym$ and $z$ 
-- which is contained in an apartment, 
as is the convex hull of any two faces -- 
is a right-angled spherical triangle 
with the combinatorial structure as depicted in the figure. 
Let $\tau$ denote the type $123$ face contained in it with vertex $z$. 
We denote the 4-vertex on $mz$ adjacent to $m$ by $w$. 
Then $\angle_m(w,x_1)=\pihalf$, 
and the arc connecting $\ora{mw}$ and $\ora{mx_1}$ in $\Si_mB$ 
consists only of a single 1-simplex of type $41$. 
Accordingly, 
$\tau':=CH(m,w,x_1)$ is a face of type $124$. 
Furthermore, 
$CH(\tau\cup\tau')=CH(x_1,y,z)$ 
is a spherical triangle 
with the combinatorial structure as shown in the figure. 
(That the geodesic triangle in $K$ with vertices $x_1,y$ and $z$ is rigid, 
follows more directly from triangle comparison,
because a comparison triangle in the unit sphere $S^2$ 
with vertices $\hat x_1,\hat y$ and $\hat z$ 
and with the same side lengths 
has angle 
$\angle_{\hat y}(\hat x_1,\hat z)=\angle_y(x_1,z)$.)
We note that $\Si_wK$ contains 
a singular 1-sphere $\kappa$ of type $\dots1212\dots$. 
Since $K$ contains no $4I$-vertices, 
this implies in particular that $dim(K)=3$. 
Then $\kappa$ bounds a 2-dimensional hemisphere in $\Si_wK$
and $w$ is a $4H$-vertex. 
In particular, we see that 
$K$ contains $4H$-vertices if it contains $1P$-vertices. 

Choosing $P:=A$, 
we infer that $K$ contains $4H$-vertices, 
and in particular, $4A$-vertices.  
The existence of $4A$-vertices in $K$ implies, dually, 
that $K$ also contains $1H$-vertices. 

Now we choose $P:=H$. 
It follows that there exists a configuration 
of $1H$-vertices $x_1,x_2\in K$ and a 1-vertex $z\in K$ 
as considered above. 
Let $u\in K$ be a 4-vertex adjacent to $x_1$ 
such that $\Si_{x_1}K$ contains a two-dimensional hemisphere 
centered at $\ora{x_1u}$.
As noted before, 
$\Si_{x_1}K\subseteq\ol B_{\pihalf}(\ora{x_1u})$. 
The spherical building $\Si_{x_1}B$ 
has Dynkin diagram 
\includegraphics[scale=0.6]{F4link1dynk.eps} and,
due to the geometry of the $B_3$-Coxeter complex, 
every 2-vertex in $\ol B_{\pihalf}(\ora{x_1u})$
is adjacent to $\ora{x_1u}$. 
Since 
$\angle_{x_1}(y,z)=\angle_{x_1}(y,u)+\angle_{x_1}(u,z)$, 
the direction $\ora{x_1u}$ must bisect $\angle_{x_1}(y,z)$ 
and thus $u=w$. 
It follows that the segment $x_1w=x_1u$ 
can be extended beyond $w=u$ inside $K$. 
Since $\ora{x_1u}$ is an interior vertex of $\Si_{x_1}K$, 
this in turn implies that $w$ is a $4I$-vertex, 
a contradiction.
\qed

\subsection{The $E_6$-case}
\label{sec:e6}

We will use the information in section~\ref{sec:e6geom}
regarding the geometry of the $E_6$-Coxeter complex. 

Let $B$ be a spherical building of type $E_6$ 
and let $K\subset B$ be a convex subcomplex. 
We denote 
$G:=Stab_{Aut(B)}(K)$ and $H:=Stab_{Inn(B)}(K)$. 
Then $H$ is a normal subgroup of $G$ and, 
in view of 
$G/H\embed Isom(\De_{mod}^{E_6})\cong\Z_2$, 
it has index $\leq2$. 
The automorphisms in $G-H$ (if any) 
preserve the vertex types 1 and 4,  
and switch the types
$2\leftrightarrow6$ and $3\leftrightarrow5$. 
We assume that $K$ is 
a counterexample to the Center Conjecture 
in the sense of section~\ref{sec:counter},
i.e.\ $K$ is no subbuilding of $B$ 
and the action $G\acts K$ has no fixed point. 

If $P$ is any property defined for $i$-vertices in $K$
and invariant under $H$,
we call vertices with this property {\em $iP$-vertices}.
Let again $A$ denote the property of 
{\em not having antipodes in $K$}. 
\begin{lem}
\label{lem:far2p}
Let $P$ be an $H$-invariant property 
defined for 2- and 6-vertices in $K$ 
and implying $A$, $P\Rightarrow A$. 
Then for any $2P$-vertex ($6P$-vertex) $x\in K$
exists another $2P$-vertex ($6P$-vertex) $x'\in Hx$
with $d(x,x')=\2pithird$. 
\end{lem}
\proof
Since the roles of 2- and 6-vertices in $E_6$-geometry are dual,
it suffices to treat the case of $2P$-vertices. 

Let us assume the contrary. 
Then the orbit $Hx$ consists of $2P$-vertices 
with pairwise distances $\arccos\quart$,
any two of which are connected by a type $232$ singular segment. 
In particular, 
$\diam(Hx)<\pihalf$ and 
$H$ has a fixed point in $K$,
cf.\ section~\ref{sec:rad}. 
Hence $G\supsetneq H$ 
and there exists $\al\in G-H$.
The orbit $H\al x$ consists of $6P$-vertices. 
Since $P\Rightarrow A$, 
none of them is antipodal to $x$.
On the other hand,
they cannot all be adjacent to $x$, 
because then $\diam(Gx)\leq\arccos\quart<\pihalf$ 
and $G$ would fix a point in $K$.
Thus there exists $y\in H\al x$ 
such that $xy$ is a type $216$ singular segment 
of length $d(x,y)=\arccos(-\quart)$. 

\parpic{\includegraphics[scale=0.4]{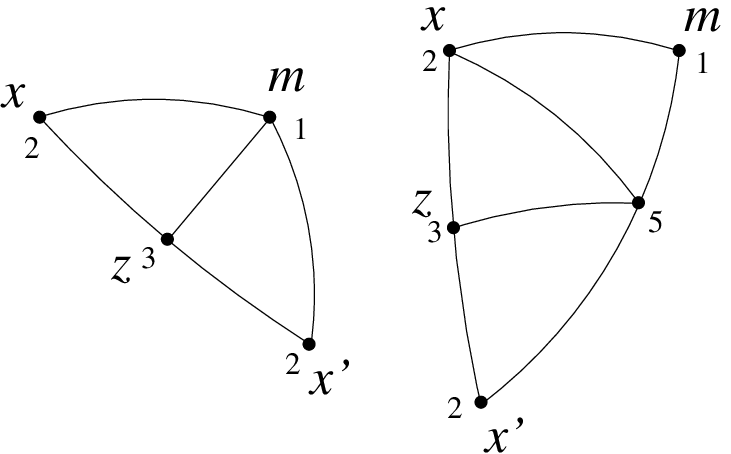}}
We claim that for the 1-vertex $m$ on $xy$ holds 
$\rad(Gx,m)\leq\pihalf$. 
By duality,
it suffices to check that $d(m,x')\leq\pihalf$ 
for all $x\neq x'\in Hx$.
We first observe that $\angle_x(m,x')<\pi$
because $d(x',y)<\pi$. 
The building $\Si_xB$ is of type $D_5$ with Dynkin diagram 
\hpic{\includegraphics[scale=0.4]{E6link2dynk_sym.eps}}.
Hence the 1-vertex $\ora{xm}$ and the 3-vertex $\ora{xx'}$
are either adjacent or connected by a type $153$ segment
in $\Si_xK$. 
In the first case,
$m$ is also adjacent to $x'$ 
and hence $d(m,x')<\pihalf$. 
(If $z$ denotes the midpoint of $xx'$, a 3-vertex,
then the building $\Si_zB$ has type $A_1\circ A_4$
with Dynkin diagram 
\includegraphics[scale=0.4]{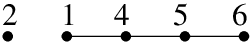}
and splits as a spherical join.)
In the other case,
the side $mx'$ of the spherical triangle $CH(m,x,x')$ 
is a singular segment of type $152$
and length $\pihalf$. 
(The link of $B$ of a 5-vertex splits off a type $A_4$-factor 
with Dynkin diagram 
\includegraphics[scale=0.4]{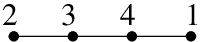}, 
and a type $2321$ segment therein has length $\pi$.)
This shows that $\rad(Gx,m)\leq\pihalf$. 

If $\rad(Gx)<\pihalf$ 
then $G$ fixes a point in $K$,
which is a contradiction.
Hence $\rad(Gx)=\pihalf$ and $m\in Cent(Gx)\cap CH(Gx)$. 
Since 1-vertices have root type,
the convex ball $\ol B_{\pihalf}(m)$ is a subcomplex 
and contains the {\em simplicial} convex hull $SCH(Gx)$ of $Gx$. 
Therefore also $\rad(SCH(Gx))=\pihalf$ and 
$m\in Cent(SCH(Gx))\cap SCH(Gx)$.
Applying Corollary~\ref{cor:fixptoncompl} to $SCH(Gx)$ 
yields that the action $G\acts K$ has a fixed point,
contradiction. 
\qed

\medskip
The midpoint of a pair of $2A$-vertices ($6A$-vertices) 
$x,x'\in K$ with distance $\2pithird$ 
is a $6A$-vertex ($2A$-vertex),
cf.\ Lemma~\ref{lem:antipodes}. 
We define the properties $M_i$ for 2- and 6-vertices in $K$ 
inductively as follows. 
Let $M_0:=A$. 
We say that a 6-vertex (2-vertex) in $K$ 
has property $M_i$, $i\geq1$, 
if it is the
{\em midpoint of a pair of $2M_{i-1}$-vertices 
($6M_{i-1}$-vertices) in $K$ with distance $\2pithird$}. 
There are the implications 
$A=M_0\Leftarrow\dots\Leftarrow M_i\Leftarrow M_{i+1}
\Leftarrow \dots$. 
We conclude from Lemma~\ref{lem:far2p}:
\begin{cor}
\label{cor:ifathenmi}
If $K$ contains a $2A$- or a $6A$-vertex,
then $K$ contains $2M_i$- and $6M_i$-vertices for all $i\geq0$.
\end{cor}

Our strategy will be 
to investigate the links of $M_i$-vertices 
for increasing $i$ 
and look for larger and larger spheres 
until we find apartments. 
We begin with the $M_2$-vertices. 
\begin{lem}
\label{lem:linkm2}
The link $\Si_wK$ of a $2M_2$-vertex $w\in K$ contains 
a type $656565656$ singular 1-sphere.
\end{lem}
\proof
Consider the following configuration. 
\parpic{\includegraphics[scale=0.3]{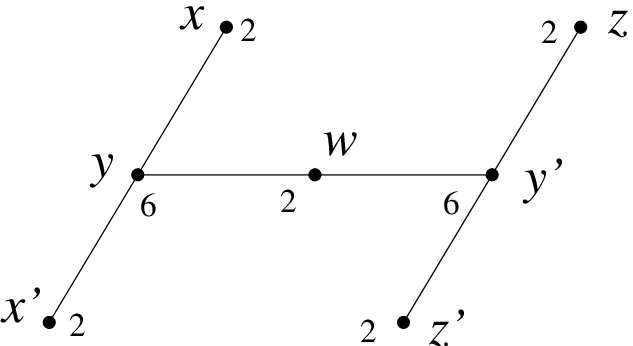}}
Let the $2M_2$-vertex $w\in K$ be the midpoint of a pair 
of $6M_1$-vertices $y,y'\in K$, 
and suppose that $y$ ($y'$) 
is the midpoint of a pair of $2A$-vertices 
$x,x'\in K$ ($z,z'\in K$). 

Since $d(x,y')<\pi$, 
we have $\angle_y(x,w)<\pi$. 
The building $\Si_yB$ is of type $D_5$ 
with Dynkin diagram
\hpic{\includegraphics[scale=0.4]{E6link6dynk_sym.eps}},
and any two distinct non-antipodal 2-vertices in it have 
distance $\pihalf$. 
Hence $\angle_y(w,x)\leq\pihalf$ and $\angle_y(w,x')\leq\pihalf$.
Since $\angle_y(x,x')=\pi$, 
we have equality 
$\angle_y(w,x)=\angle_y(w,x')=\pihalf$ 
and $\ora{yw}$ is the midpoint 
of a type $23232$ geodesic segment in $\Si_yK$ 
connecting $\ora{yx}$ and $\ora{yx'}$. 

Because of $d(w,x)<\2pithird$, 
the segment $wx$ has type $232$,
and analogously the segments $wx'$, $wz$ and $wz'$. 
\parpic{\includegraphics[scale=0.6]{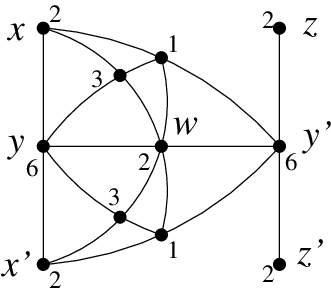}}
The 3-vertices $\ora{wx},\ora{wx'}$ 
and the antipodal 6-vertices $\ora{wy},\ora{wy'}$ 
lie on a type $6316136$ singular 1-sphere in $\Si_wK$. 
We see that the link $\Si_{yw}K$ of the edge $yw$ 
contains a pair of antipodal 3-vertices, 
while $\Si_{y'w}K$ contains a pair of antipodal 1-vertices. 
Exchanging the roles of $y$ and $y'$, 
we obtain also antipodal 1-vertices in $\Si_{yw}K$
and antipodal 3-vertices in $\Si_{y'w}K$.

We wish to produce in, say, $\Si_{yw}K$
from the pairs of antipodal 1- and 3-vertices 
a pair of antipodal 5-vertices.
Note that 
the spherical buildings $\Si_{yw}B$ and $\Si_{y'w}B$ 
have type $D_4$ and Dynkin diagram
\hpic{\includegraphics[scale=0.3]{E6link26dynk.eps}}.

\begin{sublem}
\label{sublem:findcircled4}
Let $L$ be a convex subcomplex 
of a spherical building $B'$ 
of type $D_4$ with Dynkin diagram 
\hpic{\includegraphics[scale=0.3]{E6link26dynk.eps}}.
Suppose that $L$ contains a pair of antipodal 1-vertices 
and a pair of antipodal 3-vertices. 
Then it contains a singular 1-sphere of type $1351351$. 
\end{sublem}
\proof
We denote the two antipodal 1-vertices by $a$ and $b$,
and the two antipodal 3-vertices by $c$ and $d$.
If $c$ (or $d$) lies on a minimizing segment $\ga$ of type $1351$ 
connecting $a$ and $b$, 
then a singular 1-sphere of the desired type 
is obtained as the convex hull of $d$ ($c$)
and a neighborhood of $c$ ($d$) in $\ga$. 

Suppose now that $d(a,c)+d(c,b)>\pi$. 
We first describe the convex hull of $a,b$ and $c$.
The segments $ac$ and $bc$ are of type $153$. 
We denote the 5-vertices on them by $p''$, respectively, $p'$. 
Let us consider the segments 
$ap''c''b$ and $ac'p'b$ of length $\pi$. 
The 5-vertices $p''$ and $p'$ cannot be antipodal,
because they are adjacent to $c$. 
So the 5-vertex $\ora{ap''}$ and the 3-vertex $\ora{ac'}$ 
in $\Si_aL$ 
are not antipodal and hence adjacent. 
\parpic{\includegraphics[scale=0.4]{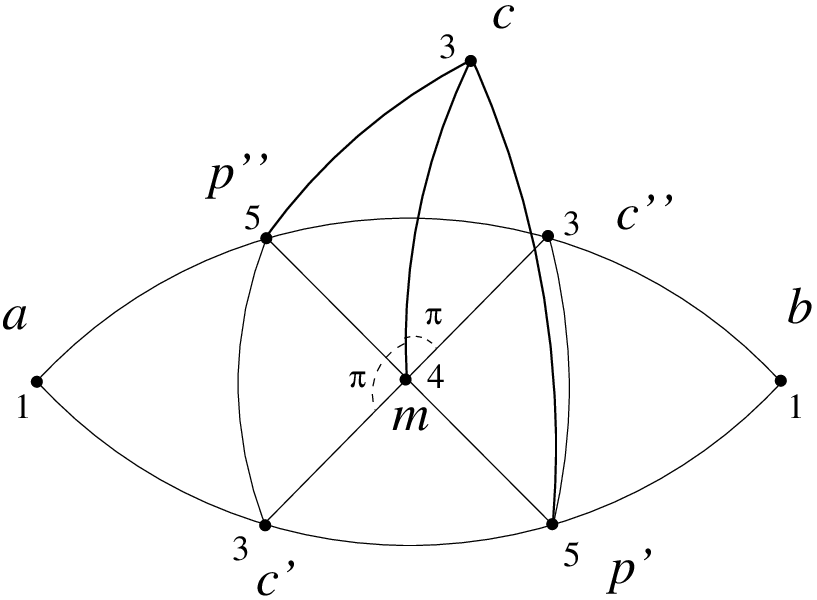}}
The convex hull of $ap''c''b$ and $ac'p'b$
is then a spherical bigon $\beta$.
Note that $\Si_{p''}B'$ is a building of type $A_3$
with Dynkin diagram
\includegraphics[scale=0.4]{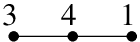}.
It contains a segment of length $\pi$ and type $1343$
from $\ora{p''a}$ to $\ora{p''c''}$ through $\ora{p''c'}$.
Its 4-vertex corresponds to a 4-vertex $m$
adjacent to $p'',c'$ and $c''$.
Replacing $p''$ by $c''$,
we see in the same way that $m$ is also adjacent to $p'$.
Hence the bigon $\beta$ is centered at $m$.
Since $c$ is adjacent to $p'$ and $p''$,
it is also adjacent to their midpoint $m$. 
The 3-vertices $\ora{mc}$, $\ora{mc'}$ and $\ora{mc''}$
are pairwise antipodal in $\Si_mL$. 
(Note that $\Si_mB'$ is a building of type $A_1\circ A_1\circ A_1$.) 

The direction $\ora{md}$ is antipodal to $\ora{mc}$
and hence antipodal to at least one of the directions 
$\ora{mc'}$ and $\ora{mc''}$. 
It follows that $d$ is antipodal to $c'$ or $c''$ 
and, as in the beginning of the proof, 
that $L$ contains the desired singular 1-sphere. 
\qed

\begin{rem}
\label{rem:findcircled4}
The proof shows that we can choose the circle in $L$ to contain the two
antipodal 1-vertices or the two antipodal 3-vertices.
(If $L$ contains a type $1351351$ circle 
through one of the two antipodal 1-vertices,
then it contains another such circle through both of them.) 
\end{rem}

\no
{\em End of proof of Lemma~\ref{lem:linkm2}.} 
Applying Sublemma~\ref{sublem:findcircled4} 
to $L=\Si_{yw}K\cong\Si_{\ora{wy}}\Si_wK$,  
we find a type $1351351$ singular 1-sphere. 
We will only use the pair of antipodal 5-vertices on it. 
These can be regarded in $\Si_wK$ 
as the directions 
$\ora{\eta\xi}$ and $\ora{\eta\xi'}$
to 5-vertices $\xi$ and $\xi'$ 
adjacent to $\eta:=\ora{wy}$.
Then $CH(\eta,\xi,\xi',\ora{wy'})$ 
is the desired type $656565656$ singular 1-sphere in $\Si_wK$. 
\qed
\begin{rem}
The 1-sphere in $\Si_{yw}K$
provided by Sublemma~\ref{sublem:findcircled4}
yields in fact a singular 2-sphere in $\Si_wK$
which contains the type $656565656$ singular 1-sphere. 
\end{rem}
Now we turn our attention to the $M_3$-vertices. 
\begin{lem}
\label{lem:linkm3}
The link $\Si_vK$ of a $6M_3$-vertex $v\in K$ contains 
a singular 2-sphere 
isomorphic to the $B_3$-Coxeter complex 
with Dynkin diagram
\includegraphics[scale=0.5]{B3_432_dynkin.eps}. 
\begin{center}
\includegraphics[scale=0.4]{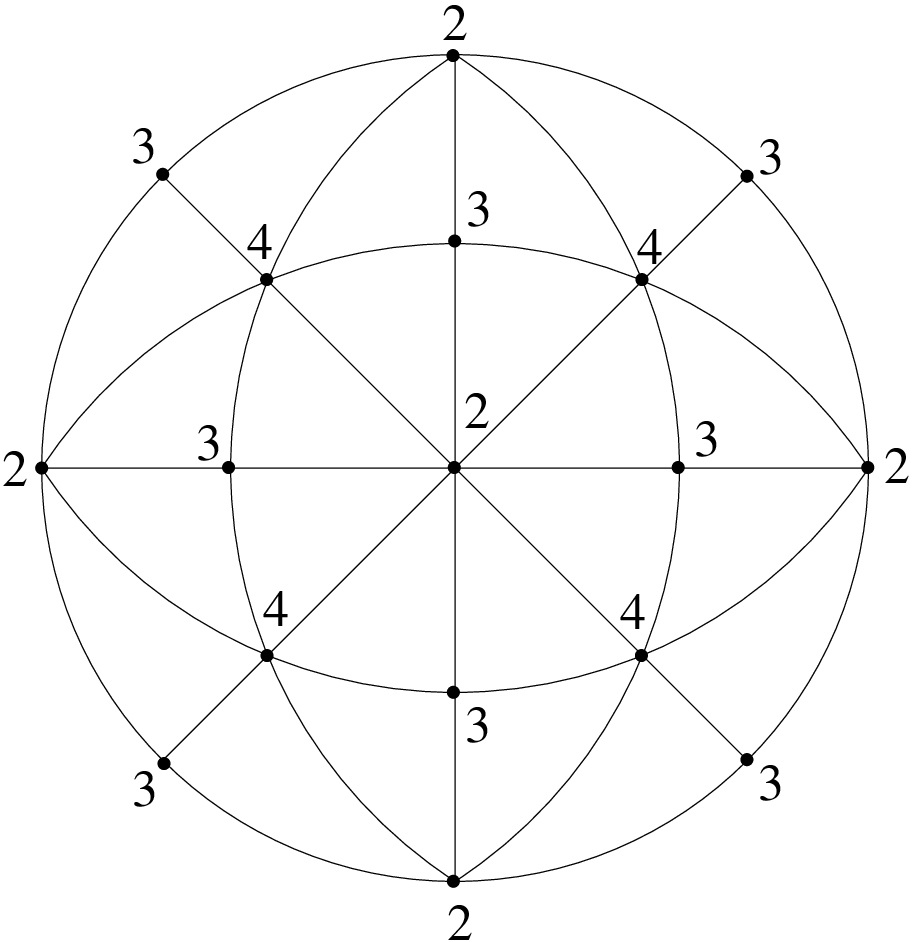}
\end{center}
\end{lem}
\proof
Let $v\in K$ be a $6M_3$-vertex
which is the midpoint of the pair of $2M_2$-vertices
$w,w'\in K$. 
Let $\ga\subset\Si_wK$ be a singular 1-sphere 
as provided by Lemma~\ref{lem:linkm2}. 
The building $\Si_wB$ is of type $D_5$ 
with Dynkin diagram
\hpic{\includegraphics[scale=0.4]{E6link2dynk_sym.eps}}.
The 6-vertices on $\ga$ are not antipodal to $\ora{wv}$,
because $w'$ has no antipodes in $K$. 
Hence they have distance $\leq\pihalf$ from $\ora{wv}$.
Since they come in (two) pairs of antipodes,
equality holds,
i.e.\ the 6-vertices on $\ga$ 
have distance $=\pihalf$ from $\ora{wv}$,
compare the beginning of the proof of Lemma~\ref{lem:linkm2}.
It follows that 
$d(\ora{wv},\cdot)\equiv\pihalf$ on $\ga$, 
and $CH(\ga\cup\{\ora{wv}\})=:h$ 
is a 2-dimensional hemisphere 
with center $\ora{wv}$ and boundary circle $\ga$. 
It is a simplicial subcomplex of $\Si_wK$ 
composed of twenty four $654$-triangles
and isomorphic to a hemisphere 
in the type $B_3$ Coxeter complex with Dynkin diagram
\includegraphics[scale=0.4]{B3_456_dynkin.eps}
as described in section~\ref{sec:e6geom}.

$\Si_{\ora{wv}}h=:\de$ is a type $545454545$ singular 1-sphere 
in $\Si_{\ora{wv}}\Si_wK\cong\Si_{wv}K\cong\Si_{\ora{vw}}\Si_vK$.
(The building $\Si_{wv}B$ has type $D_4$ with Dynkin diagram 
\hpic{\includegraphics[scale=0.3]{E6link26dynk.eps}}.) 
When regarded as a circle in the latter space,
$\de$ is the link at $\ora{vw}$ of a singular 2-sphere in $\Si_vK$
with poles $\ora{vw}$ and $\ora{vw'}$.
To determine the link $\de'$ of the 2-sphere 
at the opposite pole $\ora{vw'}$,
we recall that in the $D_5$-Coxeter complex with Dynkin diagram 
\hpic{\includegraphics[scale=0.3]{E6link6dynk_sym.eps}}
the possible types for singular segments connecting 2-antipodes 
are $24342$, $2512$ and $23232$.
Hence to a 4-vertex (5-vertex) on $\de$ 
corresponds a 4-vertex (1-vertex) on $\de'$,
i.e.\ $\de'$ has type $141414141$.
(This can also be expressed by saying that 
the natural isomorphism 
$\Si_{wv}K\cong\Si_{\ora{vw}}\Si_vK
\cong\Si_{\ora{vw'}}\Si_vK\cong\Si_{w'v}K$
of type $D_4$ spherical buildings 
with Dynkin diagram 
\hpic{\includegraphics[scale=0.3]{E6link26dynk.eps}}
switches the vertex types $1\leftrightarrow5$.)
We regard $\de'$ as a singular 1-sphere in 
$\Si_{w'v}K\cong\Si_{\ora{vw'}}\Si_vK$.

By exchanging the roles of $w$ and $w'$, 
we obtain likewise a type $141414141$ singular 1-sphere 
in $\Si_{wv}K$,
besides the type $545454545$ singular 1-sphere
which we found before. 
From these, 
we wish to produce a type $343434343$ singular 1-sphere. 
This leads us to the following continuation 
of Sublemma~\ref{sublem:findcircled4}. 
\begin{sublem}
\label{sublem:find2sphd4}
Let $L$ be a convex subcomplex 
of a spherical building $B'$ of type $D_4$ 
with Dynkin diagram 
\hpic{\includegraphics[scale=0.3]{E6link26dynk.eps}}
Suppose that $L$ contains 
a type $141414141$ singular 1-sphere 
and a pair of antipodal 5-vertices. 
Then it contains a singular 2-sphere. 
\end{sublem}
We recall that in $D_4$-geometry 
there is only one type of singular 2-spheres, 
see section~\ref{sec:e6geom},
and hence a singular 2-sphere 
contains singular 1-spheres of all possible types. 
\proof
We denote the type $141414141$ singular 1-sphere by $\eps$.
Let $g,\hat g \in \eps$ be two antipodal 1-vertices.
By Sublemma~\ref{sublem:findcircled4} and Remark~\ref{rem:findcircled4}
there is a singular 1-sphere $\eps'$ of type 1351351 containing
$g,\hat g$. (The roles of 1-,3- and 5-vertices in $D_4$-geometry
are equivalent,
as reflected by the symmetries of the Dynkin diagram.)

Since $g$ has an antipode in $L$,
it suffices to find a singular 1-sphere in $\Si_gL$.
The spherical building $\Si_gB'$
has type $A_3$
and Dynkin diagram
\includegraphics[scale=0.5]{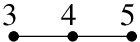}.
We know already that $\Si_gL$ contains
the pair of antipodal 4-vertices
$\{\xi,\hat\xi\}=\Si_g\eps$
and the pair of antipodes $\{\eta,\hat\eta\}=\Si_g\eps'$
consisting of a 3-vertex $\eta$ and a 5-vertex $\hat\eta$.

To find the singular 1-sphere in $\Si_gL$, 
we proceed as in the proof of Sublemma~\ref{sublem:findcircled4}. 
If $\eta$ or $\hat\eta$ lies on a minimizing segment 
connecting $\xi$ and $\hat\xi$,
then we are done. 
Otherwise, 
let us consider the convex hull of $\xi,\hat\xi$ and, say, $\eta$.
The segments $\xi\eta$ and $\hat\xi\eta$
are of type $453$. 
We call the 5-vertices on them $\zeta''$, respectively, $\zeta'$. 
They are distinct, 
and we denote by $\mu$ 
the midpoint of the type $545$ segment $\zeta''\zeta'$. 
The segments $\xi\zeta''\hat\xi$ and $\xi\zeta'\hat\xi$
have types $4534$ and $4354$, 
and the spherical bigon 
bounded by them is right angled.
(Note that 
the spherical building $\Si_{\xi}\Si_gB'$ has type $A_1\circ A_1$ 
with Dynkin diagram
\includegraphics[scale=0.5]{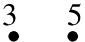}. 
It is a complete bipartite graph with edge lengths $\pihalf$.)
The 4-vertex $\mu$ is the center of the bigon.
\parpic{\includegraphics[scale=0.4]{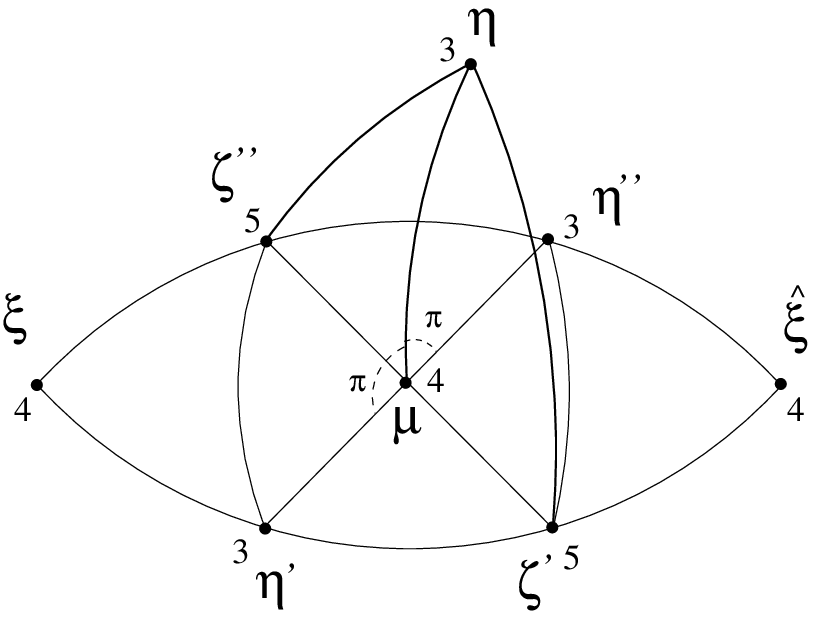}}
Since $\eta$ is adjacent to $\zeta''$ and $\zeta'$,
it is also adjacent to their midpoint $\mu$. 
The segment $\eta\mu\hat\eta$ has type $3435$. 
Let $\eta''$ and $\eta'$ 
denote the 3-vertices on the type $435$ segments 
$\hat\xi\zeta''$, respectively, $\xi\zeta'$. 
The 3-vertices $\ora{\mu\eta''}$ and $\ora{\mu\eta'}$ 
in $\Si_{\mu}\Si_gL$ are antipodal. 
The 3-vertex $\ora{\mu\hat\eta}$
is antipodal to at least one of them, say to $\ora{\mu\eta'}$.
It follows that $\eta'$ is an antipode of $\hat\eta$ 
and $CH(\hat\eta,\xi,\eta',\zeta')$
is the singular 1-sphere in $\Si_gL$ 
which we are looking for. 
\qed

\no
{\em End of proof of Lemma~\ref{lem:linkm3}.} 
Applying Sublemma~\ref{sublem:find2sphd4}
to $L=\Si_{vw}K\cong\Si_{\ora{vw}}\Si_vK$, 
we find a singular 2-sphere. 
We will only use one of the type $343434343$ singular 1-spheres 
contained in it. 
This circle can be regarded 
as the link at $\ora{vw}$ 
of a singular 2-sphere in $\Si_vK$ 
with poles the 2-vertices $\ora{vw}$ and $\ora{vw'}$.
It is a 2-sphere of the desired type. 
\qed

\begin{rem}
The 2-sphere in $\Si_{vw}K$
provided by Sublemma~\ref{sublem:find2sphd4}
yields in fact a singular 3-sphere in $\Si_vK$
which contains this singular 2-sphere. 
\end{rem}
Finally, we look at $M_4$-vertices.
\begin{lem}
\label{lem:linkm4}
The link $\Si_uK$ of a $2M_4$-vertex $u\in K$ 
contains an apartment 
and is therefore a top-dimensional subbuilding of $\Si_uB$. 
\end{lem}
\proof
Let $u\in K$ be a $2M_4$-vertex
which is the midpoint of the pair of $6M_3$-vertices
$v,v'\in K$. 
Let $s\subset\Si_vK$ 
be a singular 2-sphere 
as provided by Lemma~\ref{lem:linkm3}.
The building $\Si_vB$ is of type $D_5$ 
with Dynkin diagram
\hpic{\includegraphics[scale=0.4]{E6link6dynk_sym.eps}}.
The 2-vertices on $s$ are not antipodal to $\ora{vu}$,
because $v'$ has no antipodes in $K$. 
Hence they have distance $\leq\pihalf$ from $\ora{vu}$.
Since they come in pairs of antipodes,
equality holds,
i.e.\ the 2-vertices on $s$ 
have distance $=\pihalf$ from $\ora{vu}$,
compare the beginning of the proofs
of Lemmas~\ref{lem:linkm2} and ~\ref{lem:linkm3}.
It follows that 
$d(\ora{vu},\cdot)\equiv\pihalf$ on $s$, 
and $CH(s\cup\{\ora{vu}\})=:h'$ 
is a 3-dimensional hemisphere 
with center $\ora{vu}$ and boundary 2-sphere $s$. 

However, 
unlike the circle and 2-hemisphere found before,
$h'\subset\Si_vK$ 
is {\em not} a simplicial subcomplex.
This can be seen, for instance, 
from the fact that it contains quadruples of 2-vertices 
with pairwise distances $\pihalf$, 
cf.\ the discussion of $\Si_2$ in section~\ref{sec:e6geom}.
Accordingly, 
the 2-sphere 
$\Si_{\ora{vu}}h'\subset\Si_{\ora{vu}}\Si_vK\cong\Si_{uv}K$ 
is not a subcomplex of the 
type $D_4$ building $\Si_{uv}B$ with Dynkin diagram 
\hpic{\includegraphics[scale=0.3]{E6link26dynk.eps}}.
(It contains triples of 3-vertices 
with pairwise distances $\pihalf$.)
The simplicial convex hull of $\Si_{\ora{vu}}h'$
is an apartment 
contained in $\Si_{uv}K$.
Since $u$ is an interior point of the edge $vuv'$ 
contained in $K$,
it follows that there are apartments also in $\Si_uK$.
\qed

\medskip
Let again $I$ denote the property of 
{\em being an interior point of $K$}. 
Clearly, $I\Rightarrow A$ 
because by assumption $K$ is no subbuilding. 
Lemma~\ref{lem:linkm4} says that $M_4\Rightarrow I$
for 2- and 6-vertices. 
\begin{lem}
\label{lem:noint26}
There are no $2I$- and $6I$-vertices in $K$.
\end{lem}
\proof
Otherwise,
we suppose without loss of generality 
that $K$ contains $2I$-vertices.
By Lemma~\ref{lem:far2p}, 
there exist two $2I$-vertices in $K$ 
with distance $\2pithird$. 
Then among the 6-vertices in $K$ adjacent to one of them 
are antipodes of the other,
a contradiction. 
\qed

\begin{cor}
\label{cor:26antip}
All 2- and 6-vertices in $K$ have antipodes in $K$. 
\end{cor}
\proof
Otherwise,
Corollary~\ref{cor:ifathenmi} and 
Lemma~\ref{lem:linkm4} imply 
that $K$ contains $2I$-vertices,
which contradicts Lemma~\ref{lem:noint26}. 
\qed

\medskip
Now the main work is done 
and we enter the endgame of our argument. 
\begin{lem}
\label{lem:1antip}
All 1-vertices in $K$ have antipodes in $K$. 
\end{lem}
\proof
Suppose that $K$ contains $1A$-vertices. 
We argue as in case 2 
in the proof of Theorem~\ref{thm:tccf4}. 
The set of all $1A$-vertices in $K$, 
which have distance $\leq\pihalf$ from any other $1A$-vertex in $K$,
must be empty in view of Corollary~\ref{cor:diampihalffixpt},
because it is $G$-invariant and has diameter $\leq\pihalf$. 
Hence, there exists a pair of $1A$-vertices $x,x'\in K$ 
with distance $>\pihalf$, i.e.\ with distance $\2pithird$. 
The midpoint $y$ is another $1A$-vertex, 
and there exists yet another $1A$-vertex $z\in K$ 
with $d(y,z)=\2pithird$. 

Since $z$ is no antipode of $x$ or $x'$,
we have 
$0<\angle_y(x,z),\angle_y(x',z)<\pi$.
At least one of these angles is $\geq\pihalf$,
say, the first one. 
The building $\Si_yB$ has type $A_5$ with Dynkin diagram
\includegraphics[scale=0.4]{E6link1dynk.eps}.
The 4-vertices $\ora{yx},\ora{yz}\in\Si_yK$ 
must have distance $\arccos(-\third)$
and their simplicial convex hull is a rhombus
whose other diagonal is a $26$-edge,
compare the discussion of $\Si_1$ in section~\ref{sec:e6geom}. 
In particular,
there exists a 2-vertex $w\in K$ adjacent to $y$.
By Corollary~\ref{cor:26antip}, 
it has an antipode $\hat w\in K$. 
The segment $wy\hat w$ is of type $21656$.
Since the 6-vertex on it adjacent to $y$ 
has an antipode in $K$, too, 
it follows that $y$ has an antipode in $K$,
a contradiction. 
\qed

\begin{lem}
\label{lem:35antip}
All 3- and 5-vertices in $K$ have antipodes in $K$. 
\end{lem}
\proof
By duality, it is enough to treat the case of 3-vertices.

Suppose that $x\in K$ is a $3A$-vertex.
If $K$ contains a 1-vertex $y$ adjacent to $x$,
then it also contains an antipode $\hat y$ of $y$
(Lemma~\ref{lem:1antip}).
The segment $yx\hat y$ contains a 6-vertex $z$ adjacent to $x$.
It has an antipode $\hat z$ in $K$ 
and, by Lemma~\ref{lem:antipodes},
$x$ has an antipode in $K$, a contradiction.
The same reasoning shows that $K$ cannot contain 
6- or 2-vertices adjacent to $x$,
i.e.\ it contains at most 4- and 5-vertices adjacent to $x$. 

\parpic(7cm,0pt)[l]{\includegraphics[scale=0.5]{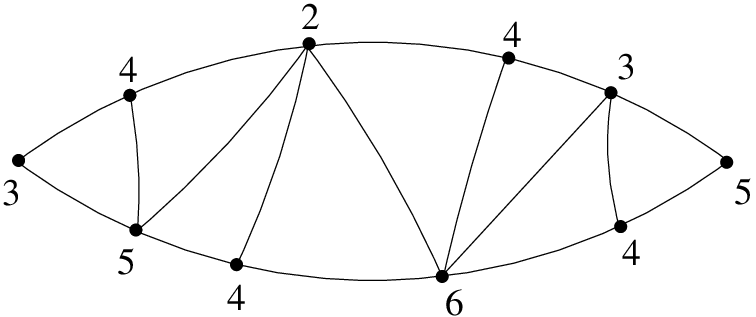}}
Hence, 
for any point $p\in K$ 
the direction $\ora{xp}$ lies on a $45$-edge in $\Si_xB$.
As a consequence, 
$p$ is contained in a spherical bigon $\beta\subset B$ 
with $x$ as one of its tips and with $\Si_x\beta$ a $45$-edge. 
(Of course, $\beta\not\subset K$.)
It has the combinatorial structure as shown in the figure. 
(This is easily verified by taking into account 
that the links of 4-vertices in $B$ have type $A_2\circ A_1\circ A_2$
with Dynkin diagram 
\includegraphics[scale=0.4]{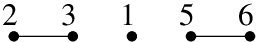}.)
If $K$ contains another 3-vertex $x'$, 
then $xx'$ can only be a type $34243$ singular segment. 
Since the 2-vertex on it has an antipode in $K$
(Corollary~\ref{cor:26antip}),
it follows that also $x$ has an antipode in $K$,
a contradiction. 

Thus $x$ is the only 3-vertex in $K$
and must be fixed by $H$. 
Since $G\acts K$ does not have a fixed point,
we have $G\supsetneq H$ and $Gx$ consists of $x$ and a $5A$-point.
They cannot be antipodal
and their unique midpoint is fixed by $G$,
a contradiction. 
\qed

\begin{lem}
\label{lem:4antip}
All 4-vertices in $K$ have antipodes in $K$. 
\end{lem}
\proof
Suppose that $x\in K$ is a $4A$-vertex.
Then all vertices in $K$ adjacent to $x$ 
have antipodes in $K$. 
If $K$ contains vertices adjacent to $x$, 
then the same reasoning as in the beginning of the proof 
of Lemma~\ref{lem:35antip}
shows that $x$ has an antipode, too, 
a contradiction. 
Hence $dim(K)=0$, 
which is also a  contradiction. 
\qed

\medskip
We proved that all vertices in $K$ have antipodes in $K$.
With Proposition~\ref{lem:allvertantip} 
it follows that $K$ is a subbuilding. 
This contradicts our assumption that $K$ is a counterexample 
to the Center Conjecture,
and we obtain our second main result: 
\begin{thm}
\label{thm:tcce6}
The Center Conjecture~\ref{conj:tcc} 
holds for spherical buildings of type $E_6$. 
\end{thm}

\subsection{The case of classical types}
\label{sec:class}

The Center Conjecture for the spherical buildings 
of classical types ($A_n$, $B_n$ and $D_n$)
was first proven by M\"uhlherr and Tits in \cite{MuehlherrTits} using 
combinatorial methods and the incidence geometries of the respective buildings.
We present in this section a proof from the point of view of CAT(1) spaces
and using methods of comparison geometry,
compare \cite[ch.\ 4.1]{diss}.

L.\ Kramer informed us about a similar argument in the $A_n$-case
showing that a convex subcomplex $K\subset B$ 
is a subbuilding or $Stab_{Inn(B)}(K)$ fixes a point in $K$. 

We will use the information in section~\ref{app:coxeter}
regarding the geometry of the Coxeter complexes. 

Let $K$ be a convex subcomplex of a spherical building $B$ of classical type
and suppose that it is not a subbuilding.
By Proposition~\ref{lem:allvertantip}, 
there are vertices in $K$ 
without antipodes in $K$. Let 
$t=\max\,\{\;i\;|\;\exists \; iA\text{-vertex in }K\}$ 
(As before,
we call an $i$-vertex of $K$ without antipodes in $K$ an {\em $iA$-vertex}.)

\begin{lem}\label{lem:minimal}
Let $B$ be of type $A_n$ or $B_n$ and
let $x\in K$ be a $tA$-vertex.
Then
there is no vertex of type $>t$ in $K$ adjacent to $x$ 
\end{lem}
\proof
Let $t'>t$ and suppose that there exists a $t'$-vertex $y\in K$ adjacent to $x$.
The maximality of $t$ implies that $y$ has an antipode $\widehat y\in K$.
The edge $yx$ extends to the singular segment
$yx\widehat y\subset K$ of length $\pi$.
Notice that $\Si_x B$ splits off a factor of type $A_{n-t}$ 
and its
Dynkin diagram has labels $t+1,\dots,n$. This implies that the direction
$\ora{x\widehat y}$ has type $t''>t$. 
It follows that the segment $x\widehat y\subset K$ has a $t''$-vertex $z$ 
in its interior and, 
by Lemma~\ref{lem:antipodes}, $z$ cannot have antipodes in $K$,
contradicting the maximality of $t$. 
\qed

\subsubsection{The $A_n$-case}

\begin{thm}
The Center Conjecture~\ref{conj:tcc} 
holds for spherical buildings of type $A_n$.
\end{thm}
\proof
We assume that $n\geq2$,
because otherwise the assertion is trivial. 

Let $K$ be a convex subcomplex of a spherical building $B$ of type $A_n$ 
and suppose that it is not a subbuilding.
Let 
$t_1=\min\,\{\;i\;|\;\exists \; iA\text{-vertex in }K\}$ 
and $t_2=\max\,\{\;i\;|\;\exists \; iA\text{-vertex in }K\}$.
Let $x_i\in K$ be a $t_iA$-vertex.
By Lemma~\ref{lem:minimal},
there is no vertex of type $>t_2$ in $K$ adjacent to $x_2$ 
and, analogously, no vertex of type $<t_1$ in $K$ adjacent to $x_1$. 

If $t_1=t_2$, we may choose $x_1=x_2$. 
It follows that $dim(K)=0$ and the Center Conjecture holds trivially 
in this case.
We therefore assume in the following that $t_1<t_2$. 

Consider the segment $x_1x_2$ 
as embedded in the vector space realization of the
Coxeter complex of type $A_n$ described in section~\ref{app:coxeter}, such
that $x_1=v_{t_1}$ 
(we work with vectors representing vertices)
and such that the initial part of
$x_1x_2$ is contained in the fundamental Weyl chamber $\Delta$. 
Then $x_2$ lies in the convex hull of $\De$ 
and the antipode of $x_1$ in the Coxeter complex,
a bigon which is given by all inequalities (\ref{ineq:wcan}) except ($t_1$).
The coordinates of $x_2$ are a permutation of the coordinates of $v_{t_2}$.
It follows from the observation above,
that the face of $\Delta$ spanned by
the initial part of $x_1x_2$ contains no vertices of types $1,\dots,t_1-1$.
This implies for $x_2=(a_0,\dots,a_n)\in \R^{n+1}$ that 
$a_0=\dots=a_{t_1-1}$ and $a_{t_1}\leq\dots\leq a_n$.
Note that $x_2$ is adjacent to $x_1$ 
if and only if in addition $a_{t_1-1}\leq a_{t_1}$ holds. 
If $x_2$ is not adjacent to $x_1$, 
it follows that $a_0=t_2$ and 
$$
x_2=(\underbrace{t_2,\dots,t_2}_{t_1},
	\underbrace{-(n+1-t_2),\dots,-(n+1-t_2)}_{t_2},
	\underbrace{t_2,\dots,t_2}_{n+1-t_1-t_2}).
$$
In particular,
$n+1-t_2\geq t_1$. 
Since $x_1$ and $x_2$ are not antipodal, 
we even have the strict inequality 
\begin{equation}
\label{ineq:strict}
n+1 > t_1+t_2.
\end{equation}
Consider now the embedding of $x_1x_2$ into the Coxeter complex
such that $x_2=v_{t_2}$ and
the initial part of
$x_2x_1$ is contained in $\Delta$. 
The observation above implies now that the face of $\Delta$
spanned by the initial part of 
$x_2x_1$ contains no vertices of types $t_2+1,\dots,n$.
This implies for $x_1=(b_0,\dots,b_n)\in \R^{n+1}$ that
$b_0\leq\dots\leq b_{t_2-1}$ and $b_{t_2}=\dots =b_n$. 
If $x_1$ is not adjacent to $x_2$, 
equivalently,
if $b_{t_2-1}>b_{t_2}$ 
it follows that 
$$
x_1=(\underbrace{-(n+1-t_1),\dots,-(n+1-t_1)}_{t_1+t_2-(n+1)},
	\underbrace{t_1,\dots,t_1}_{n+1-t_1},
	\underbrace{-(n+1-t_1),\dots,-(n+1-t_1)}_{n+1-t_2})
$$
and $t_1\geq n+1-t_2$.
But this inequality contradicts (\ref{ineq:strict}). 
Hence, $x_1$ and $x_2$ must be adjacent.

We saw that any $t_1A$-vertex is adjacent to any $t_2A$-vertex. 
Their distance is $\leq\diam(\De)=:\de<\pihalf$.
This implies that the set of $t_iA$-vertices 
has circumradius $\leq\de$ 
and therefore a unique circumcenter $c_i$.
Moreover, 
$c_i$ is contained in the closed convex hull 
of the set of $t_iA$-vertices,
in particular $c_i\in K$.
It follows that $c_i$ has distance $\leq\de$ 
from every $t_{3-i}A$-vertex,
and $d(c_1,c_2)\leq\de$. 

The Dynkin diagram for $A_n$ has only one nontrivial symmetry
which exchanges the labels $i\leftrightarrow (n+1-i)$. 
Hence a building automorphism in 
$Stab_{Aut(B)}(K)-Stab_{Inn(B)}(K)$ 
must switch the labels $t_1\leftrightarrow t_2$
(according to their definition)
and exchange $c_1\leftrightarrow c_2$,
whereas the automorphisms in $Stab_{Inn(B)}(K)$ fix both $c_1$ and $c_2$. 
It follows that the midpoint $m(c_1,c_2)\in K$ of $c_1$ and $c_2$ 
is fixed by the entire group $Stab_{Aut(B)}(K)$. 
\qed

\subsubsection{The cases $B_n$ and $D_n$}

\begin{thm}
The Center Conjecture~\ref{conj:tcc} 
holds for spherical buildings of type $B_n$.
\end{thm}
\proof
If $n=2$, then $\dim(K)\leq1$ 
and the Center Conjecture holds.
So, let $K$ be a convex subcomplex of a spherical building $B$ of type $B_n$ 
for $n\geq 3$
and suppose that it is not a subbuilding.
Let 
$t=\max\,\{\;i\;|\;\exists \; iA\text{-vertex in }K\}$.
Let $x\in K$ be a $tA$-vertex.
By Lemma~\ref{lem:minimal}, 
there are no vertices of type $>t$ in $K$
adjacent to $x$.

Let $x'\in K$ be another $tA$-vertex.
Consider the segment $xx'$ as embedded in the vector space realization
of the Coxeter complex of type $B_n$ described in section~\ref{app:coxeter}.
Assume that 
\begin{equation*}
x=v_t=(0,\dots,0,\underbrace{1,\dots,1}_{n+1-t})
\end{equation*} 
and that the initial part of
$xx'$ is contained in the fundamental Weyl chamber $\Delta$.
Then the coordinates of $x'$ 
satisfy all inequalities (\ref{ineq:wcbn}) except ($t$).
They agree up to permutation and signs 
with the coordinates of $x$. 
The observation above implies that 
the face of $\Delta$ spanned by the initial part of
$xx'$ contains no vertices of types $t+1,\dots,n$.
For $x'=(a_1,\dots,a_n)$ 
this means that $a_t=\dots=a_n$
(besides $0\leq a_1\leq\dots\leq a_{t-1}$ which we will not use). 
If $a_t=1$, then $x=x'$; if $a_t=0$, then $d(x,x')=\pihalf$; and if
$a_t=-1$, then $x$ and $x'$ are antipodal.
The last case cannot occur,
and hence $d(x,x')\leq\pihalf$. 
It follows that 
the $Stab_{Aut(B)}(K)$-invariant set of $tA$-vertices in $K$
has diameter $\leq\pihalf$. 
Therefore, $Stab_{Aut(B)}(K)$ fixes a point in $K$ by
Corollary~\ref{cor:diampihalffixpt}. 
\qed

\begin{thm}
The Center Conjecture~\ref{conj:tcc} 
holds for spherical buildings of type $D_n$.
\end{thm}
\proof 
For $n\geq5$, 
the $D_n$-case of the Center Conjecture follows 
from the $B_n$-case, 
because a spherical building of type $D_n$ 
can be regarded as a (thin) spherical building of type $B_n$. 
In the same vein, 
the $D_4$-case follows from the $F_4$-case,
compare (\ref{eq:autdnbn}) and (\ref{eq:autd4f4})
in section~\ref{app:coxeter}. 
\qed

\medskip
A direct proof of the $D_{n\geq5}$-case can be given as well.
We skip it here because it is very similar 
to the argument in the $B_n$-case.
To keep our treatment of the classical types self-contained,
we include a direct proof in the $D_4$-case.

\medskip
\noindent
{\em Direct proof in the $D_4$-case.}
Let $K$ be a convex subcomplex of a spherical building $B$ of type $D_4$ 
and suppose that $K$ is a counterexample to the Center Conjecture.

Suppose first, 
that $K$ contains $3A$-vertices. Recall that the 3-vertices in $D_4$
are the vertices of root type. 
Their possible mutual distances are 
$0,\pithird,\pihalf,\2pithird,\pi$. 
The midpoint of a segment connecting two 3-vertices
at distance $\pithird$ lies in the interior 
of a simplex of type 124 adjacent to both 3-vertices.

Arguing as in the beginning of case 2 in the proof of Theorem~\ref{thm:tccf4}, 
since $K$ is a counterexample, 
there exist $3A$-vertices $x,x'\in K$ at distance
$\2pithird$. 
The simplicial convex hull of the segment $xx'$ is 3-dimensional and the
midpoint $y_1$ of $xx'$ is an interior $3A$-vertex of $K$.
Let $y_2\in K$ be a $3I$-vertex at distance $\2pithird$ to $y_1$.
Since $y_i$ is interior, 
we can find $z_i\in K$, with $d(z_i,y_i)=\frac{\pi}{6}$,
such that $\ora{y_iz_i}$ is antipodal to $\ora{y_iy_{3-i}}$ 
in $\Si_{y_i}K$ for $i=1,2$.
In particular, $z_1$ and $z_2$ are antipodal. 
Notice that $z_i$ lies in the interior 
of a simplex of type 124. It follows that $K$ contains a 2-sphere, 
contradicting Corollary~\ref{cor:cod1sph}. 
Hence all 3-vertices in $K$ have antipodes in $K$.

Since $K$ is a counterexample, 
there is a vertex $w\in K$ without antipodes in $K$,
cf.\ Proposition~\ref{lem:allvertantip}. 
It has type $i\in\{1,2,4\}$. 
If there exists a 3-vertex $v\in K$ adjacent to $w$,
then it has an antipode $\hat v$ in $K$. 
The interior of the segment $w\hat v$ 
contains another 3-vertex $v'$. 
By Lemma~\ref{lem:antipodes}, 
$v'$ is a $3A$-vertex, contradiction. 
Thus $w$ cannot be adjacent to a 3-vertex in $K$.
This implies that $w$ is the only $i$-vertex in $K$,
because any two distinct nonantipodal 
$i$-vertices are connected by a segment of type $i3i$.

It follows that the non-empty $Stab_{Aut(B)}(K)$-invariant set $V$ 
of $A$-vertices in $K$
consists of vertices of pairwise different types $\neq3$.
In particular, $|V|\leq3$
and the possible distances between vertices in $V$ 
are $\pithird$ and $\2pithird$. 
If $\diam(V)>\pihalf$,
we may assume without loss of generality
that $V$ contains a $1A$-vertex $u_1$ and a $2A$-vertex $u_2$
with distance $\2pithird$. 
The segment $u_1u_2$ then has type $142$.
Its midpoint is a $4A$-vertex $u_4$
and $V=\{u_1,u_2,u_4\}$. 
We conclude that always $\rad(V)<\pihalf$.
Hence $V$ has a unique circumcenter in $K$
and it is fixed by $Stab_{Aut(B)}(K)$. 
\qed

\subsection{An intrinsic version of the results}

Let us briefly indicate that 
in all cases considered in this paper 
our arguments actually yield a 
slightly stronger intrinsic version of the Center Conjecture:
\begin{thm}\label{thm:generalversion}
Suppose that $B$ is a spherical building of type $F_4$, $E_6$
or of classical type 
and that $K\subseteq B$ is a convex subcomplex. 
Then $K$ is a subbuilding 
or the action $Aut(K)\acts K$ 
of the automorphisms of $K$ has a fixed point. 
\end{thm}

By an {\em automorphism} of a convex subcomplex $K\subset B$ we mean an 
isometry of $K$ preserving the polyhedral
structure induced by $B$ 
with the additional restriction
that the change of vertex types under the automorphism 
corresponds to a symmetry of the Dynkin diagram of $B$.
Note that an automorphism of a subcomplex 
does in general not extend to an
automorphism of the ambient building. 

To see that our arguments 
yield this more precise version of the results,
one first observes that Corollaries~\ref{cor:diampihalffixpt},
\ref{cor:fixptoncompl} and \ref{cor:cod1sph} 
hold by the same reasons for the larger group action 
$Aut(K)\acts K$ 
(see also Remark~\ref{rem:fixedpt}). 
One modifies the notion of {\em counterexample} 
as introduced in section~\ref{sec:counter} 
by replacing the action $Stab_{Aut(B)}(K)\acts K$ with $Aut(K)\acts K$ 
and it follows that counterexamples keep the properties 
discussed there. 
Let $Inn(K)\subset Aut(K)$ denote 
the subgroup of automorphisms of $K$ 
which preserve vertex types. 
Then, after replacing $Stab_{Aut(B)}(K)$ with $Aut(K)$
and $Stab_{Inn(B)}(K)$ with $Inn(K)$, 
the proofs in sections~\ref{sec:f4}, \ref{sec:e6} and \ref{sec:class}
yield Theorem~\ref{thm:generalversion}.

\end{document}